\documentclass[twoside,11pt]{article}
\usepackage{amssymb}
\usepackage[dvips]{graphics}
\newtheorem{definition}{Definition}
\newtheorem{theorem}{Theorem}
\newtheorem{conjecture}{Conjecture}
\newtheorem{claim}{Claim}
\newtheorem{lemma}{Lemma}
\def\C{\mathbb{C}}
\def\R{\mathbb{R}}
\def\Z{\mathbb{Z}}

\def\Q{\mathbb{Q}}

\def\ra{{\longrightarrow}}
\def\la{{\longleftarrow}}
\def\c{\mathcal{C}}
\def\om{{\otimes\dots\otimes}}
\def\g{\mathfrak{g}}
\def\H{\mathcal{H}}
\def\G{{\Gamma}}
\def\eps{{\epsilon}}
\def\O{\mathcal{O}}
\def\E{\mathcal{E}}
\def\F{\mathcal{F}}
\def\tc{{\tilde {C}}}
\def\A{\mathbb{A}}

\begin{document}
\title{{\huge{\bf{Operads and Motives \\in Deformation Quantization}}} 
\footnote
{To be published in Letters in Mathematical Physics {\bf 48} (1), April 1999}}
\author{\sc{\large{ Maxim Kontsevich}}\\
{\it Institut des Hautes \'Etudes Scientifiques}\\
{\it 35, route de Chartres, F-91440 Bures-sur-Yvette, France}\\
 {\small{e-mail:}} {\tt {\small{maxim@ihes.fr}}}}
\date{April 13, 1999\\ \vskip5mm
\centerline{\it Dedicated to the memory of Mosh\'e Flato}}
\maketitle
\thispagestyle{empty}
\begin{abstract}
The algebraic world of associative algebras has many deep connections with the 
geometric world of two-dimensional surfaces. Recently D.~Tamarkin discovered 
that the operad of chains of the little discs operad is formal, i.e.
it is homotopy equivalent to its cohomology. From this fact and from  
Deligne's conjecture on Hochschild complexes follows almost immediately my 
formality result in deformation quantization. I review the situation as it 
looks now. Also I conjecture that the motivic Galois group acts on deformation 
quantizations, and speculate on possible relations of higher-dimensional 
algebras and of motives to quantum field theories.
\end{abstract}

\section{Introduction}

Deformation quantization was proposed about 22 years ago in the 
pioneering work of Bayen, Flato, Fr\o nsdal, Lichnerowicz, and 
Sternheimer [BFFLS] as an alternative to the usual correspondence 
\begin{eqnarray*}
\mathrm{classical\ systems}&\leftrightarrow &\mathrm{quantum\
  systems} \\
\mathrm{symplectic\ manifolds}&\leftrightarrow &\mathrm{Hilbert\ spaces}
\end{eqnarray*}

The idea is that algebras of observables in quantum mechanics
are ``close to" commutative algebras of functions on manifolds (phase spaces). 
In other words, quantum algebras of observables are \textit{deformations}
of commutative algebras.
      
In the first order in perturbation theory one obtains automatically
a Poisson structure on phase space. Remind that a Poisson structure 
on a smooth manifold $X$ is a bilinear operation $\{\cdot,\cdot\}$
on $C^{\infty}(X)$ satisfying the Jacobi identity
and the Leibniz rule with respect to the usual product in $C^{\infty}(X)$.
Two typical examples of Poisson manifolds are symplectic manifolds
and dual spaces to Lie algebras. A star-product on a Poisson manifold $X$
is an associative (but possibly non-commutative) product on $C^{\infty}(X)$ 
depending formally on a parameter usually denoted by~$\hbar$ 
(the ``Planck constant''). The product should have the form
\begin{equation}
f\star g= fg+\hbar\{f,g\}+\sum_{n\ge 2}\hbar^n B_n(f,g)
\end{equation}
where $(B_n)_{n\ge 1}$ are bidifferential operators
$C^{\infty}(X)\otimes C^{\infty}(X)\ra C^{\infty}(X)$.
On the set of star-products acts an infinite-dimensional gauge group
of linear transformations of the vector space $C^{\infty}(X)$
depending formally on~$\hbar$, of the form
\begin{equation}
f\mapsto f +\sum_{n\ge 1} \hbar^n D_n(f)
\end{equation}
where the $(D_n)_{n\ge 1}$ are differential operators on $X$.
 
About two years ago I proved (see [K1]) that for every 
Poisson manifold $X$ there is a canonically defined
gauge equivalence class of star-products on $X$.
This result  gave a complete answer to first basic
questions in  the program started by Mosh\'e Flato and co-authors.
   
I obtained the existence of a canonical deformation quantization
from  a more general and stronger result, the formality theorem.
The statement of this theorem is that in a suitably defined
homotopy category of differential graded Lie algebras, 
two objects are equivalent. The first object is the Hochschild complex
of the algebra of functions on the smooth manifold $X$, and the
second object is a $\Z$-graded Lie superalgebra of polyvector fields
on $X$. In the course of the proof I constructed an explicit isomorphism
in the homotopy category of Lie algebras for the case
$X=\R^n$. The terms in the formula for this isomorphism
can be naturally identified with Feynman diagrams
for certain two-dimensional quantum field theory
with broken rotational symmetry (see [CF] for a detailed
derivation of Feynman rules for this theory).
The moral is that a kind of ``string theory'' is necessary
for deformation quantization, which was originally
associated with quantum mechanics! 
Mosh\'e was extremely happy when I told him about
my construction and we celebrated the  ending of an old story
with a bottle of champagne in his Paris
apartment. Mosh\'e  was  very enthusiastic 
about the new approach and expected new developments.
              
Later it became clear that not only  there exists a
canonical way to quantize, but that one can define a universal
infinite-dimensional manifold parametrizing quantizations.
There are several evidences that this universal
manifold is a principal homogeneous space of the 
so-called Grothen\-dieck-Teichm\"uller group, introduced by Drinfeld and 
Ihara. At the ICM98 Congress in Berlin, I gave a talk 
(which Mosh\'e qualified as ``wild") on relations between deformations, 
motives and the Grothendieck-Teichm\"uller group. After the Congress 
I decided not to write notes of the talk because one month later a 
dramatic breakthrough in the area happened, which shed light on some parts 
of my talk and required  further work; the present paper fills this gap. 

Dmitry Tamarkin found a new short derivation of the formality 
theorem (for the case $X=\R^n$) from a very general result concerning 
all associative algebras. Also from his result the (conjectured) 
relation between deformation quantizations and the 
Grothendieck-Teichm\"uller group  seems to be much more
transparent. The present paper is a result of  my attempts to 
understand and generalize Tamarkin's results.
I use all the time the language of operads
and of homotopy theory of algebraic structures.
For me it seems to be the first real application of operads
to algebraic questions.

The main steps of the new proof of the formality theorem in
deformation quantization are the following:

\noindent STEP 1) On the cohomological Hochschild complex of any 
associative algebra acts the operad of chains in the little discs operad 
(Deligne's conjecture). This result is purely topological/combinatorial.
  
\noindent STEP 2) 
 The operad of chains of the little discs operad is formal, i.e.
 it is quasi-isomorphic to its cohomology. More precisely, this is true
 only in characteristic zero, e.g. over the field $\Q$ of rational numbers.
  All the proofs of existence of a quasi-isomorphism 
  use certain multi-dimensional integrals and give explicit 
  formulas over real or over complex numbers.
 
\noindent STEP 3) From steps 1) and 2) follows that for any algebra $A$
over a field of characteristic zero, its cohomological
Hochschild complex $C^*(A,A)$ and its Hochschild cohomology $B:=H^*(A,A)$
are algebras over the \textit{same} operad (up to homotopy). 
Moreover, if one chooses an explicit homotopy of complexes
between the Hochschild complex of $A$ and the graded space $B$ considered
as a complex with zero differential, one obtains \textit{two}
different structures of a homotopy Gerstenhaber algebra on $B$.
If these two structures are not equivalent, the first non-zero
obstruction to the equivalence gives a non-zero element in the
$H^1(\mathsf{Def}(B))$ of the deformation complex of the 
Gerstenhaber algebra $B$.
    
\noindent STEP 4) For the case $A=\R[x_1,\dots,x_n]$, the algebra $B$
is the Gerstenhaber algebra of polynomial polyvector fields on $\R^n$. 
An easy calculation shows that $H^1(\mathsf{Def}(B))$ coincides with
$B^2$ (i.e. with the space of bivector fields on $\R^n$).
The explicit homotopy between  $B$ and $C^*(A,A)$ can be made invariant
under the group $\mathsf{Aff}(\R^n)$ of affine transformations of $\R^n$.
There is no non-zero $\mathsf{Aff}(\R^n)$-invariant bivector field on $\R^n$.
The conclusion is that two structures of the homotopy Gerstenhaber
algebra on $B$, mentioned in step 3), are equivalent.
In particular, these two structures give equivalent structures
of homotopy Lie algebra on $B$. This implies that
$C^*(A,A)$ is  equivalent to $B$ as a homotopy Lie algebra,
which is the statement of the formality result in [K1].
    
Step 3) is the main discovery of D.~Tamarkin. It is an absolutely fundamental
 new fact about \textit{all} associative algebras.
 It implies in particular that the Hochschild complex is quasi-isomorphic
  to another natural complex with strictly associative and commutative
   cup-product. This other complex is hard to write down explicitly 
at the moment. Also, there are many different
 ways to identify operads as in step 2), even  up to homotopy.
 All these choices form an infinite-dimensional algebraic manifold defined
 over the field of rational numbers $\Q$. On this manifold acts 
 the Grothendieck-Teichm\"uller group.
 By the naturalness of the calculation in step 4), the same group
 acts on quantizations of Poisson structures on $\R^n$, and also
 on the set of gauge equivalence classes of star-products in $\R^n$.
Presumably, the action extends to general Poisson manifolds and to
general star-products.

The paper is organized as follows:

\noindent Section 2: an introduction to  operads, and to
 Deligne's conjecture concerning the Hochschild complex.
  I would like to apologize for certain vagueness in subsections 2.5-2.7.

\noindent Section 3: a proof of the result of Tamarkin on the
formality of the chain operad of small discs operad, and a sketch of
its application to my formality theorem in deformation quantization.
In fact the idea of the proof presented here goes back to 1992-1993,
but somehow at  that time I  missed the point.

\noindent Section 4: I describe in elementary terms a version of
motives and of the motivic Galois group, and indicate its relations with
homogeneous spaces appearing in various questions in deformation
quantization.

\noindent Section 5: speculations about the possible r\^ole in Quantum
Field Theories of things described in previous sections.
 
\section{Deligne's conjecture and its generalization \\ to higher 
         dimensions}
\subsection{Operads and algebras}    
Here I remind the definitions of an operad and of an algebra over an
 operad (see also [GJ], [GK]).  The language of operads is  convenient
 for descriptions and constructions of various algebraic structures.
 It became quite popular in theoretical physics during the past few years
 because of the emergence of many new types of algebras related with
 quantum field theories.

\begin{definition}
An operad (of vector spaces) consists of the following:

\noindent 1) a collection of vector spaces $P(n)$, $n\ge 0$,

\noindent 2) an action of the symmetric group $S_n$ on $P(n)$ for every $n$,

\noindent 3) an identity element $\mathrm{id}_P\in P(1)$,

\noindent 4) compositions $m_{(n_1,\dots, n_k)}$:
\begin{equation}
P(k)\otimes \bigl(P(n_1)\otimes P(n_2)\otimes\dots\otimes 
P(n_k)\bigr) \ra P(n_1+\dots+n_k)
\end{equation} 
for every $k\ge 0$ and $n_1,\dots,n_k\ge 0$
satisfying a natural list of axioms which will be clear from examples. 
\end{definition}

The simplest example of operad is given by 
$P(n):= \mathsf{Hom}(V^{\otimes n}, V)$ where $V$ is a vector space. 
The action of the symmetric group and the identity element are obvious, 
and the compositions are defined by the  substitutions
 $$\bigl(m_{(n_1,\dots, n_k)}\bigl(\phi\otimes (\psi_1\otimes \psi_2\otimes
\dots \otimes \psi_k)\bigr)
 \bigr)(v_1\otimes\dots\otimes v_{n_1+\dots +n_k})$$
 $$:= \phi\bigl(\psi_1(v_1\otimes\dots \otimes v_{n_1})\otimes\dots
 \otimes \psi_k(v_{n_1+\dots+n_{k-1}+1}\otimes\dots
\otimes v_{n_1+\dots + n_k})\bigr)$$
where $\phi\in P(k)=\mathsf{Hom}(V^{\otimes k},V),$ $\psi_i\in
P(n_i)=\mathsf{Hom}(V^{\otimes n_i},V)$, $i=1,\dots,k$. 
 
 This operad is called the endomorphism operad of a vector space. The
axioms in the definition of operads express natural properties of this
example. Namely, there should be an associativity law for multiple
compositions, various compatibilities for actions of symmetric groups, and
evident relations for compositions including the identity element.
 
Another important example of an operad is $\mathsf{Assoc}_1$. 
The $n$-th component $\mathsf{Assoc}_1(n)$ for $n\ge 0$ is defined as 
the collection of all universal ( = functorial) $n$-linear operations
$A^{\otimes n}\ra A$ defined on all associative algebras $A$ with unit.
The space $\mathsf{Assoc}_1(n)$ has dimension $n!$, and is spanned
by the operations
 $$a_1\otimes a_2\otimes\dots\otimes a_n\mapsto a_{\sigma(1)}a_{\sigma(2)}
\dots a_{\sigma(n)}$$
where $\sigma\in S_n$ is a permutation. The space $\mathsf{Assoc}_1(n)$ 
can be identified with the subspace of the free associative unital algebra 
in $n$ generators consisting of expressions polylinear in each generator.
        
\begin{definition}
An algebra over an operad $P$ consist of a vector space
 $A$ and a collection of polylinear maps $f_n:P(n)\otimes A^{\otimes n}\ra A$
  for all $n\ge 0$ satisfying the following list of axioms:
  
\noindent  1) for any $n\ge 0$  the map  $f_n$ is $S_n$-equivariant,
  
\noindent   2) for any $a\in A$ we have $f_1(\mathrm{id}_P\otimes a)=a$,
  
\noindent   3) all compositions in $P$ map to compositions
   of polylinear operations on $A$. 
\end{definition}
   
In other words, the structure of algebra over $P$  on a vector space $A$
is given by a homomorphism of operads from $P$ to the endomorphism
operad of $A$. Another name for algebras over $P$ is $P$-algebras.
   
For example, an algebra over the operad $\mathsf{Assoc}_1$ is an associative 
unital algebra. If we replace the $1$-dimensional space $\mathsf{Assoc}_1(0)$ 
by the zero space $0$, we obtain an operad $\mathsf{Assoc}$ describing 
associative algebras possibly  without unit. Analogously, there is an 
operad denoted $\mathsf{Lie}$, such that  $\mathsf{Lie}$-algebras  are 
Lie algebras. The dimension of the $n$-th component $\mathsf{Lie}(n)$ is
$(n-1)!$ for $n\ge 1$ and $0$ for $n=0$.
 
 Let us warn the reader that not all algebraic structures correspond
 to operads. Two examples of classes of algebraic structures that cannot
be cast in the language of operads are the class of fields and the 
class of Hopf algebras.
 
 We conclude this section with an explicit description of free algebras in
terms of operads.
 
\begin{theorem}  
Let $P$ be an operad and $V$ be a vector space. Then the free
$P$-algebra $\mathsf{Free}_P(V)$ generated by $V$ is naturally isomorphic as a
vector space to 
$$\bigoplus_{n\ge 0} \bigl(P(n)\otimes V^{\otimes n}\bigr)_{S_n}$$
\end{theorem}
  
The subscript $\cdot_{S_n}$ denotes the quotient space of coinvariants for 
the diagonal action of group $S_n$. The free algebra $\mathsf{Free}_P(V)$  
is defined by the usual categorical adjunction  property:
 
%\parbox[t]{110mm}
{\it the set $\mathsf{Hom}_{\,P-\mathrm{algebras}}
(\mathsf{Free}_P(V),A)$ (i.e. homomorphisms in the category of $P$-algebras)

is naturally isomorphic to the set 
$\mathsf{Hom}_{\,\mathrm{vector\,\,spaces}} (V,A)$ for any $P$-algebra $A$.}

\subsection{Topological operads,  operads of complexes, etc.}    

In the definition of operads we can replace vector spaces by 
topological spaces, and the operation of tensor product by the usual
(Cartesian) product.
 
\begin{definition} 
A topological operad consists of the following:
 
\noindent 1) a collection of topological spaces $P(n)$, $n\ge 0$,

\noindent 2) a continuous action of the symmetric group 
$S_n$ on $P(n)$ for every $n$,

\noindent 3) an identity element $\mathrm{id}_P\in P(1)$,

\noindent 4) compositions $m_{(n_1,\dots, n_k)}$:

$P(k)\times \bigl(P(n_1)\times P(n_2)\times\dots\times 
P(n_k)\bigr) \ra P(n_1+\dots+n_k)$  \\
which are continuous maps for every $k\ge 0$ and $n_1,\dots ,n_k\ge 0$
satisfying a list of axioms analogous to the one in the definition of an
operad of vector spaces.
\end{definition}
 
An analog of the endomorphism operad is the following one:
for any $n\ge 0$ the topological space $P(n)$ is the space of continuous 
maps from $X^n$ to $X$, where $X$ is a given compact topological space.
  
In general, the definitions of an operad and of an algebra over an operad 
can be made in arbitrary symmetric monoi\-dal category $\c$
(i.e in a category endowed with the functor $\otimes:\c\times \c\ra\c$, 
the identity element $1_{\c}\in \mathsf{Objects}(\c)$, and various coherence 
isomorphisms for associativity, commutativity of $\otimes$, etc.,  see, 
for example,  [ML]).
     
We shall here consider mainly operads in the symmetric monoidal
category $\mathsf{Complexes}$ of $\Z$-graded complexes of abelian groups 
(or vector spaces over a given field). Often operads in the category of 
complexes are called \textit{differential graded operads}, or simply 
\textit{dg-operads}. Each component $P(n)$ of an operad of complexes is a 
complex, i.e. a vector space decomposed into a direct sum
$P(n)=\oplus_{i\in \Z} P(n)^i$, and endowed with a differential
$d:P(n)^i\ra P(n)^{i+1}$, $d^2=0$, of degree $+1$.
     
There is a natural way to construct an operad of complexes from a
topological operad by using a version of the singular chain complex.
Namely, for a topological space $X$, denote by $\mathsf{Chains}(X)$ 
the complex concentrated in nonpositive degrees, whose $(-k)$-th component 
for $k=0,1,\dots$ consists of the formal finite additive combinations
$$\sum_{i=1}^N n_i \cdot f_i,\,\,\,n_i\in \Z,\,\,\,N\in \Z_{\ge 0}$$
of  continuous maps $f_i:[0,1]^k\ra X$ (``singular cubes'' in $X$), 
modulo the following relations 
      
\noindent 1) $f\circ \sigma=\mathit{sign}(\sigma) \, f$ for
any $\sigma \in S_k$ acting on the standard cube
$[0,1]^k$ by permutations of coordinates, 
      
\noindent 2) $f'\circ pr_{k\ra (k-1)}=0$ where 
$pr_{k\ra (k-1)}:[0,1]^k\ra[0,1]^{k-1}$ is the projection onto first 
$(k-1)$ coordinates, and $f':[0,1]^{k-1}\ra X$ is a continuous map.
       
The boundary operator on cubical chains is defined in the usual way.
The main advantage of cubical chains with respect to simplicial 
chains is that there is an external product map
\begin{equation}
\bigotimes_{i\in I} \bigl(\mathsf{Chains}(X_i)\bigr)\ra 
\mathsf{Chains}\bigl(\prod_{i\in I} X_i\bigr)
\end{equation}
which is a \textit{natural} homomorphism of complexes  for any finite 
collection $(X_i)_{i\in I} $ of topological spaces.
        
Now if $P$ is a topological operad then the collection of complexes 
$\bigl(\mathsf{Chains}(P(n))\bigr)_{n\ge 0}$ has a natural operad
structure in the category of complexes of Abelian groups. The compositions
in $\mathsf{Chains}(P)$ are defined using the external tensor product 
of cubical chains.

 Passing from complexes to their cohomology we obtain an operad 
$H_*(P)$ of $\Z$-graded Abelian groups ( = complexes with zero differential), 
the homology operad of $P$.

\subsection{Operad of little discs}    

Let $d\ge 1$ be an integer. Denote by $G_d$ the $(d+1)$-dimensional Lie
group acting on $\R^d$ by affine transformations $u\mapsto \lambda u +v$ 
where $\lambda>0$ is a real number and $v\in \R^d$ is a vector. 
This group acts simply transitively on the space of closed discs in
$\R^d$ (in the usual Euclidean metric). The disc with center $v$ and 
with radius $\lambda$ is obtained from the standard disc 
$$D_0:=\{(x_1,\dots,x_d)\in \R^d | \,x_1^2+\dots +x_d^2\le 1\}$$
by a transformation from $G_d$ with parameters $(\lambda,v)$.
 
\begin{definition}
The little discs operad $C_d$ is a topological operad with the 
following structure:
 
\noindent  1) $C_d(0)=\emptyset$,
 
\noindent 2) $C_d(1)= \,\mathrm{point}\,=\{\mathrm{id}_{C_d}\}$,
 
\noindent 3) for $n\ge 2$ the space $C_d(n)$ is the space of configurations 
of $n$ disjoint discs $(D_i)_{1\le i\le n}$ inside the standard disc $D_0$.
  
The composition  $C_d(k)\times C_d(n_1)\times\dots\times C_d(n_k)\ra
C_d(n_1+\dots+n_k)$ is obtained by applying elements from 
 $G_d$ associated with  discs $(D_i)_{1\le i\le k}$ in the configuration 
in $C_d(k)$ to configurations in all $C_d(n_i),\,\,i=1,\dots,k$
and putting the resulting configurations together. The action of 
the symmetric group $S_n$ on $C_d(n)$ is given by renumerations of 
indices of discs $(D_i)_{1\le i\le n}$.
\end{definition}
  
The operad $C_d$ was introduced in 70-ies by Boardmann and Vogt, 
and by Peter May (see [BV], [M]) in order to describe homotopy types of 
$d$-fold loop spaces (i.e. spaces of continuous maps
$$\mathsf{Maps}((S^d,{\rm base \,\,\, point}), (X,x))$$ 
where $X$ is a topological space with base point $x$). It is the most 
important operad in homotopy theory. Strictly speaking, topologists use a 
slightly different model called the operad of little cubes, but the 
difference is not essential because homotopically  there is no difference 
between cubes and discs.
   
The space $C_d(n)$ is homotopy equivalent to the configuration space
of $n$ pairwise distinct points in $\R^d$:
\begin{displaymath}
\mathsf{Conf}_n(\R^d):=\quad \left(\R^d\right)^n\setminus \mathrm{Diag}=
\{(v_1,\dots,v_n)\in \left(\R^d\right)^n |v_i\ne v_j 
\ \mathrm{for\ any}\ i\ne j\}
\end{displaymath}
    
There is an obvious continuous map $C_d(n)\ra \mathsf{Conf}_n(\R^d)$ which 
associates  a collection of disjoint discs with the collection of
their centers. This map induces a homotopy equivalence because each fiber
of this map is contractible.

The space $\mathsf{Conf}_2(\R^d)$ (and hence $C_d(2)$) is homotopy
equivalent to the $(d-1)$-dimensional sphere $S^{d-1}$.
The homotopy equivalence is given by the map
$$(v_1,v_2)\mapsto {v_1-v_2\over |v_1 -v_2|}\in S^{d-1}\subset \R^d$$
        
\subsection{Hochschild complex and Deligne's conjecture }    
 
In 1993 Pierre Deligne made a conjecture relating the little discs operad 
(in dimension $d=2$) and the cohomological Hochschild complex
of an arbitrary associative algebra $A$ (defined over any field $k$). 
Namely, the Hochschild complex $C^*(A,A)$ is concentrated in non-negative 
degrees and is defined as
\begin{equation}
C^n(A,A):=\mathsf{Hom}_{\mathrm{vector\,spaces}} (A^{\otimes n}, A),
\,\,\,n\ge 0
\end{equation}
and the differential in $C^*(A,A)$ is given by the formula
\begin{eqnarray*}
&&(d \phi)(a_1\otimes \dots \otimes a_{n+1}):=\\
&&a_1 \phi(a_2\otimes\dots
\otimes a_n)+\sum_{i=1}^n (-1)^i \phi(a_1\otimes
   \dots\otimes a_i a_{i+1}\otimes \dots\otimes a_{n+1})+\cdots\\
&&\quad \cdots +(-1)^{n+1} \phi(a_1\otimes \dots\otimes a_n) a_{n+1}
\end{eqnarray*}
   for any $\phi\in C^n(A,A)$.
   
The Hochschild complex plays a fundamental r\^ole in the deformation theory 
of an associative algebra $A$. There are two basic operations on 
$C^*:=C^*(A,A)$, the cup product $\cup:C^k\otimes C^l\ra C^{k+l}$  
and the Gerstenhaber bracket $[\,,\,]:C^k\otimes C^l\ra C^{k+l-1}$. 
The formulas for these operations are the following (here $\phi\in C^k$, 
$\psi\in C^l$):
\begin{eqnarray*}
&&(\phi\cup\psi)(a_1\otimes\dots\otimes a_{k+l})\\
&&\quad :=(-1)^{kl}
\phi(a_1\om a_k)\cdot\psi(a_{k+1}\om a_{k+l})\\
&&[\phi,\psi]:=\phi\circ\psi- (-1)^{(k-1)(l-1)}
\psi\circ\phi\quad \mathrm{where}\\
&&(\phi\circ\psi)(a_1\om a_{k+l-1}):=\\
&&\quad \sum_{i=1}^{k-1}(-1)^{i(l-1)}
\phi(a_1\om a_i\otimes\psi(a_{i+1}\om a_{i+l}) \om a_{k+l-1})
\end{eqnarray*}
     
The Gerstenhaber bracket gives (after a shift of the $\Z$-grading 
by~$1$) the structure of differential graded Lie algebra on the 
Hochschild complex. The cup product has also a remarkable property: 
it is \textit{not} graded commutative (it is associative only), but 
the induced operation on the cohomology \textit{is} graded commutative.
Moreover, the  Gerstenhaber bracket induces an operation on  $H^*(A,A)$
which satisfies the Leibniz rule with respect to the cup product:
$$[\phi,\psi_1\cup \psi_2]=[\phi,\psi_1]\cup\psi_2 +(-1)^{(k-1)(l_1-1)} 
\psi_1\cup [\phi,\psi_2],\,\,\,\phi\in C^k,\,\psi_i \in C^{l_i}$$
although this identity does not hold on the level of cochains.
The cohomology space $H^*(A,A)$ carries the structure 
of the Gerstenhaber algebra, i.e. it is a graded vector space endowed 
with a Lie bracket of degree $(-1)$ (satisfying the skew symmetry 
condition and the Jacobi identity with appropriate signs),
and with a graded commutative associative product of degree $0$, 
satisfying the graded Leibniz rule with respect to the bracket.
        
It was observed by several people (F.~Cohen in [C], P.~Deligne,...)
that the $\Z$-graded operad $\mathsf{Gerst}$ describing Gerstenhaber algebras
has  a very beautiful topological meaning. Namely, it is naturally
equivalent to the homology operad of the topological operad $C_2$. The space
of binary operations $C_2(2)$ is homotopy equivalent to the circle $S^1$. The
Gerstenhaber bracket corresponds to the generator of $H_1(S^1)\simeq \Z$ 
and the cup product corresponds to the generator of $H_0(S^1)\simeq \Z$.
     
\begin{conjecture}[P.~Deligne] 
There exists a natural action of the operad $\mathsf{Chains}(C_2)$ on the 
Hochschild complex $C^*(A,A)$ for arbitrary associative algebra $A$.
\end{conjecture}

 The story of this conjecture is quite dramatic. In 1994 E.~Getzler and 
J.~Jones posted on the e-print server a preprint [GJ] in which the proof 
of the Deligne conjecture was contained.
Essentially at the same time M.~Gerstenhaber and A.~Voronov in [GV] made
analogous claims. The result was considered as well established and was 
actively used later. But in the spring of 1998 D.~Tamarkin observed that 
there was a serious flaw in both preprints. The cell decomposition used 
there turned out to be not compatible with the operad structure; the first 
example of wrong behavior  appears for operations with $6$ arguments.
       
I think that I have now a complete proof of the Deligne conjecture
(and its generalization, see the next section).
A combination of two results of Tamarkin (see [T1] and [T2])
also implies Deligne conjecture. Sasha Voronov says that he corrected
the problem in his approach, and there is also an announcement [MS]  by 
J.~McClure and J.~Smith with the same result.
       
Unfortunately, all ``proofs'' and announcements of proof are still too 
complicated to be put here. We need a really short and convincing argument 
for this very fundamental fact about Hochschild complexes.
It seems that the simplicity of the Hochschild complex is deceiving.

\subsection{Higher-dimensional generalization of  Deligne's
conjecture}    

I propose here an ``explanation'' of Deligne's conjecture and its natural
generalization.

The operad $\mathsf{Assoc}$ has a topological origin, it is the 
homology operad of the little intervals (i.e. $1$-dimensional discs)
operad. Namely, for  $n\ge 1$, the spaces $C_1(n)$ have $n!$ connected 
components corresponding to permutations $\sigma\in S_n$, and each of 
the components is contractible.
    
We can go still one step down, defining the operad $C_0$ as a really
trivial operad:
   $$C_0(n)=\emptyset\,\,\,\,{\rm for}\,\,\,n\ne
1,\,\,\,\,\,C_0(1)=\mathrm{point}=\{\mathrm{id}_{C_0}\}\,\,$$
This definition is natural because one cannot put $\ge 2$ disjoint
zero-dimensional discs inside one zero-dimensional disc ( = point).
   
\begin{definition}
For $d\ge0$, a $d$-algebra is an algebra over the operad 
$\mathsf{Chains}(C_d)$ in the category of complexes. 
\end{definition}

For $d>0$ this notion was introduced by Getzler and Jones.
By definition, a $0$-algebra is just a complex. The notion of 
$1$-algebra  is very close to the notion of an associative algebra. 
There is a well-known companion of associative algebras, a class of 
so called $A_{\infty}$-algebras. One can show that any $1$-algebra carries 
a natural structure of $A_{\infty}$-algebra, and from the points of view of
homotopy theory and of deformation theory, there is no difference
between associative algebras, $A_{\infty}$-algebras and $1$-algebras 
(see 2.6 and 3.1). In particular, one can introduce the cohomological 
Hochschild complex for $A_{\infty}$-algebras and $1$-algebras. 
This Hochschild complex carries a natural structure of a differential 
graded Lie algebra.
      
An $A_{\infty}$-version of Deligne's conjecture says that this Hochschild
complex carries naturally a structure of $2$-algebra, extending the
structure of differential graded Lie algebra. It has a  baby version 
in dimension $(0+1)$: if $A$ is a vector space (i.e. $0$-algebra 
concentrated in degree 0) then the Lie algebra of the group of
affine transformations 
$$\mathsf{Lie}(\mathsf{Aff}(A))=\mathsf{End}(A)\oplus A$$
has also a natural structure of an associative algebra, in particular it is
a $1$-algebra. The product in $\mathsf{End}(A)\oplus A$ is given by the
formula
$$(\phi_1,a_1)\times (\phi_2,a_2):=(\phi_1 \phi_2, \phi_1 (a_2))\,\,\,.$$
The space $\mathsf{End}(A)\oplus A$ plays the r\^ole of the Hochschild 
complex in the case $d=0$ (see 2.7).
         
The definition of a $d$-algebra given above seems to be too complicated 
for doing concrete calculations. At the end of section 3.2 we describe 
much smaller operads which are quasi-isomorphic to the chain operads 
of little discs operads. The reader will see that in dimension  
$d\ge 2$ the situation is very simple.
            
Now we introduce the  notion of \textit{action} of a $(d+1)$-algebra 
on a $d$-algebra. It is convenient to formulate it using so-called 
colored operads. Instead of defining exactly what a colored operad is,
we give one typical example: there is a colored operad with two colors such 
that algebras over this operad are pairs $(\g,A)$ where $\g$ is a Lie algebra
and $A$ is an associative algebra on which $\g$ acts by derivations.
         
Let us fix a dimension $d\ge 0$. Denote by $\sigma:\R^{d+1}\ra \R^{d+1}$ 
the reflection 
$$(x_1,\dots,x_{d+1})\mapsto (x_1,\dots,x_d,-x_{d+1})$$ at the coordinate
hyperplane, and by $\H_+$ the upper-half space 
         $$\{(x_1,\dots,x_{d+1})|\,x_{d+1}> 0\}$$
     
\begin{definition} 
For any pair of non-negative integers $(n,m)$ we define a topological space
$SC_d(n,m)$ as   
    
\noindent 1) the empty space $\emptyset$ if $n=m=0$,
      
\noindent 2) the one-point space if $n=0$ and $m=1$,
      
\noindent 3) in the case $n\ge 1$ or $m\ge 2$, 
the space of configurations of $m+2n$ disjoint discs 
$(D_1,\dots,D_{m+2n})$ inside the standard disc $D_0\subset \R^{d+1}$
such that $\sigma(D_i)=D_i$ for $i\le m$,  $\sigma(D_i)=D_{i+n}$ for
$m+1\le i\le m+n$ and such that all discs $D_{m+1},\dots,D_{m+n}$ are 
in the upper half space $\H_+$.
\end{definition}

The reader should think about points of $SC_d(n,m)$ as about configurations
of $m$ disjoint semidiscs $(D_1\cap \H_+,\dots,D_m\cap\H_+)$ and of
$n$ discs $(D_{m+1},\dots,D_{m+n})$ in the standard semidisc $D_0\cap\H_+$. 
The letters ``$SC$'' stand for ``Swiss Cheese'' [V]. 
Notice that the spaces $SC_d(0,m)$ are naturally isomorphic to
$C_d(m)$ for all $m$. One can define composition maps analogously to 
the case of the operad $C_d$:
$$SC_d(n,m)\times \bigl(C_{d+1}(k_1)\times\cdots\times C_{d+1}(k_n)\bigr)
\times\bigl(SC_d(a_1,b_1)\times\cdots\times SC_d(a_m,b_m)\bigr)$$
   $$\ra SC_d(k_1+\cdots+k_n+a_1+\cdots+a_m,b_1+\cdots+b_m)$$
            
\begin{definition}[A~.Voronov [V{]}] 
The colored operad $SC_d$ has two colors and consists of collections
of spaces  
$$\bigl(SC_d(n,m)\bigr)_{n,m\ge 0},\quad \bigl(C_{d+1}(n)\bigr)_{n\ge 0},$$
and appropriate actions of symmetric groups, identity elements, and
of all composition maps. 
\end{definition}

As before, we can pass from a colored operad of topological spaces to
a colored operad of complexes using the functor $\mathsf{Chains}$.   
       
\begin{definition} 
An action of a $(d+1)$-algebra $B$ on a $d$-algebra $A$ is, on the pair
$(B,A)$, a structure of algebra over the colored operad 
$\mathsf{Chains}(CS_d)$, compatible with the structures of algebras on 
$A$ and on $B$.
\end{definition}

The \textit{generalized Deligne conjecture} says that for every
$d$-algebra $A$ there exists a universal (in an appropriate sense,
up to homotopy, see 2.6) $(d+1)$-algebra acting on $A$. 
I think that I have a proof of this conjecture, so I am making the following
\begin{claim} 
Let $A$ be a $d$-algebra for $d\ge 0$.
Consider the homotopy category of pairs $(B,\rho)$ where $B$ is a 
$(d+1)$-algebra and $\rho$ is an action of $B$ on $A$.
In this category there exists a final object.
\end{claim}
  
Using this claim one can give the  
  
\begin{definition} 
For a $d$-algebra $A$ the (generalized) Hochschild complex
$\mathsf{Hoch}(A)$ is defined as the universal $(d+1)$-algebra acting on $A$.
\end{definition}

The universality of $\mathsf{Hoch}(A)$ means that in the sense 
of homotopy theory of algebras,  an action of $(d+1)$-algebra
$B$ on $A$ is the same as a homomorphism of $(d+1)$-algebras
$B\ra \mathsf{Hoch}(A)$. 
   
In the  next subsection we review homotopy and deformation theory of 
algebras over operads, and in the following one we describe an 
``explicit'' model for $\mathsf{Hoch}(A)$. 

\subsection{Homotopy theory and deformation theory}    

Let $P$ be an operad of complexes, and $f:A\ra B$ be a morphism of two
$P$-algebras. 
\begin{definition} 
$f$ is a quasi-isomorphism iff it induces an isomorphism of the cohomology 
groups of $A$ and $B$ considered just as complexes. 
\end{definition}

 Two algebras $A$ and $B$ are called \textit{homotopy equivalent}
iff there exists a chain of quasi-isomorphisms
\begin{equation}
A=A_1\ra A_2\la A_3\ra \dots\la A_{2k+1}=B
\end{equation}
  
  One can define a new structure of category on the collection of
$P$-algebras in which quasi-isomorphic algebras become equivalent.
 There are several ways to do it, using either Quillen's machinery
 of homotopical algebra (see [Q]), or using a  free resolution of 
 the operad $P$, or some simplicial constructions, etc.
 For example, in the category of differential graded Lie algebras,
 morphisms in the homotopy category are so called $L_{\infty}$-morphisms 
(see e.g. [K1]), modulo a suitably defined defined equivalence relation 
(a homotopy between  morphisms). 
 I shall not give here any precise definition of the homotopy
 category in general, just say that morphisms in the homotopy category
 of $P$-algebras are connected components of certain topological spaces, 
exactly as in the usual framework of homotopy theory (i.e. in the
category of topological spaces).

In the case when the operad $P$ satisfies some technical conditions, 
one can transfer the structure of a $P$-algebra by quasi-isomorphisms 
of complexes. In particular, one can make the construction described 
in the following lemma:
      
\begin{lemma} 
Let $P$ be an operad of complexes, such that if we consider $P$ as an operad
just of $\Z$-graded vector spaces, it is free and generated by operations 
in $\ge 2$ arguments. Let $A$ be an algebra over $P$, and let us choose
a splitting of $A$ considered as a $\Z$-graded space into the direct sum
$$A=H^*(A)\oplus V\oplus V[-1],\,\,\,(V[-1])^k:=V^{k-1}$$
endowed with a differential of the form 
$d(a\oplus b\oplus c)=0\oplus 0\oplus b[-1]$.
Then there is a canonical structure of a $P$-algebra on the cohomology 
space $H^*(A)$, and this algebra is homotopy equivalent to $A$. 
\end{lemma}           

We shall use it later (in 3.5) in combination with the fact that the 
operad $\Q\otimes \mathsf{Chains}(C_d)$ is free as an operad of $\Z$-graded
vector spaces over $\Q$. This  is evident because the action of $S_n$ on
$C_d(n)$ is free and the composition morphisms in $C_d$ are embeddings.

One can associate with any operad $P$ of complexes and, with any $P$-algebra 
$A$, some differential graded Lie algebra (or more generally, a
$L_{\infty}$-algebra, see 3.1) $\mathsf{Def(A)}$. This Lie algebra is 
defined canonically up to a quasi-isomorphism (the same as up to a homotopy).
It controls the deformations of $P$-algebra structure on $A$.
There are several equivalent constructions of $\mathsf{Def}(A)$ using
either resolutions of $A$ or resolutions of the operad $P$.
Morally, $\mathsf{Def}(A)$ is the Lie algebra of derivations in homotopy 
sense of $A$. For example, if $P$ is an operad with zero differential then
$\mathsf{Def}(A)$ is quasi-isomorphic to the differential graded Lie algebra
of derivations of $\tilde{A}$ where  $\tilde{A}$ is any free resolution of $A$. 

The differential graded Lie algebras $\mathsf{Def}(A_1)$ and 
$\mathsf{Def}(A_2)$ are quasi-isomorphic for homotopy equivalent 
$P$-algebras $A_1$ and $A_2$.
   
\subsection{Hochschild complexes and deformation theory}    
 
First of all, if $A$  a $d$-algebra then the shifted complex $A[d-1]$,
$$(A[d-1])^k:=A^{(d-1)+k}$$
carries a natural structure of $L_{\infty}$-algebra. It comes from a  
homomorphism of operads in homotopy sense from the twisted by $[d-1]$
operad $\mathsf{Chains}(C_d)$ to the operad $\mathsf{Lie}$. In order to 
construct such a homomorphism one can use  fundamental chains of all 
components of the Fulton-MacPherson operad (see 3.3.1), or deduce the 
existence of a  homomorphism from the results in 3.2.
       
Moreover, $A[d-1]$ maps as homotopy Lie algebra to $\mathsf{Def}(A)$, i.e.
$A[d-1]$ maps to ``inner derivations'' of $A$. These inner derivations 
form a Lie ideal in $\mathsf{Def}(A)$ in homotopy sense.
   
\begin{claim} 
The quotient homotopy Lie algebra $\mathsf{Def}(A)/A[d-1]$ is naturally 
quasi-isomorphic to $\mathsf{Hoch}(A)[d]$. 
\end{claim}

 In the case when $d=0$ and the complex $A$ is concentrated in degree $0$, 
 the Lie algebra $\mathsf{Def}(A)$ is $\mathsf{End}(A)$, i.e. the Lie 
algebra of linear transformations in the vector space $A$.
      
   \begin{lemma}
The Hochschild complex of $0$-algebra $A$ is $A\oplus \mathsf{End}(A)$ 
(placed in degree $0$).
   \end{lemma}
   
 Here follows a sketch of the proof. First of all, the colored operad $SC_0$ 
is quasi-isomorphic to its zero-homology operad  $H_0(SC_0)$ because all 
connected components of spaces $(SC_0(n,m))_{n,m\ge 0}$ and of 
$(C_1(n))_{n\ge 0}$ are contractible. By general philosophy (see 3.1) 
this implies that we can replace $SC_0$ by $H_0(SC_0)$ in the definition 
of the Hochschild complex given in 2.5. The $H_0$-version of a $1$-algebra 
is an associative non-unital algebra, and the $H_0$-version of an action
is the following:
           
1) an action of an associative non-unital algebra $B$ on vector space $A$
   (it comes from the generator of $\Z=H_0(SC_0(1,1)$),
             
2) a homomorphism from $B$ to $A$ of $B$-modules
(coming from the  generator of $\Z=H_0(SC_0(1,0)$).
              
It is easy to see that to define an action as above is the same as to 
define a homomorphism of non-unital associative algebras from $B$ to 
$\mathsf{End}(A)\oplus A$. Thus, the Hochschild complex is (up to homotopy) 
equal to $\mathsf{End}(A)\oplus A$.
   
 Let us continue the explanation of the Claim  2 above for the case $d=0$.  
As (homotopy) Lie algebra $\mathsf{Hoch}(A)$ coincides with the Lie algebra 
of affine transformations on $A$. The homomorphism 
$\mathsf{Def}(A)\ra \mathsf{Hoch}(A)$ is a \textit{monomorphism}, but
in homotopy category every morphism of Lie algebras can be replaced 
by an epimorphism! The Abelian Lie superalgebra $A[-1]$ is the ``kernel''
of this morphism. Let us show explicitly how all this works.
The Lie algebra $\mathsf{Def}(A)=(A)$ is quasi-isomorphic  
to the following differential graded Lie algebra $\g$: as $\Z$-graded space
it is $$\mathsf{End}(A)\oplus A\oplus A[-1],$$
i.e. the graded components of $\g$ are
$\g^0=\mathsf{End}(A)\oplus A,\,\,\,\g^1=A,\,\,\,\g^{\ne 0,1}=0$. 
The nontrivial components of the Lie bracket on $\g$ are
the usual bracket on $\mathsf{End}(A)$ and the action of $\mathsf{End}(A)$ 
on $A$ and on $A[-1]$. The only nontrivial component of the differential 
on $\g$ is the shifted by [1] identity map from $A$ to $A[-1]$.
The evident homomorphism $$\g\ra \mathsf{End}(A)$$
is a homomorphism of differential graded Lie algebras, and
also a quasi-isomorphism.
There is a short exact sequence of dg-Lie algebras
$$0\ra A[-1]\ra \g\ra \mathsf{End}(A)\oplus A\ra 0$$
          
In the case  $d=1$ an analogous thing happens.
The deformation complex of an associative algebra $A$ is
the following \textit{subcomplex} of the shifted by [1] Hochschild complex:
$$\mathsf{Def}(A)^n:=\mathsf{Hom}_{\mathrm{vector\,spaces}}
(A^{\otimes(n+1)},A)\,\,\,{\rm for}\,\,n\ge 0; \,\,\,\,
\mathsf{Def}^{<0}(A):=0$$ 
The deformation complex is quasi-isomorphic  to an $L_{\infty}$-algebra 
$\g$ which as $\Z$-graded vector space is 
 $$\mathsf{Def}(A)\oplus A\oplus A[1].$$
   The Hochschild complex of $A$ is a \textit{quotient} complex of $\g$  
by the homotopy Lie ideal $A$.
     
\section{Formality of operads of chains, application to\\
deformation quantization}    

\subsection{Quasi-isomorphisms of operads}    

Operads themselves are algebras over a certain colored operad, 
$\mathsf{Operads}$. This is quite obvious because an operad is just a 
collection of vector spaces and polylinear maps between these spaces 
satifying some identities.  If we work in characteristic zero, it is
convenient to associate colors with all Young diagrams, i.e. with all 
irreducible representations of all finite symmetric groups $S_n,\,\,n\ge 0$.
   
   Analogously to the case of algebras, we can speak about
quasi-isomorphisms of operads of complexes.

\begin{definition}
A morphism $f:P_1\ra P_2$ between two operads of complexes is called 
a quasi-isomorphism iff the maps of complexes $f(n):P_1(n)\ra P_2(n)$ 
induce isomorphisms of cohomology groups for all $n$.
\end{definition}
   
   Homotopy categories and deformation theories of algebras over 
   quasi-isomorphic operads are equivalent.
  A typical example: the operad $\mathsf{Lie}$ is quasi-isomorphic to the
   operad $L_{\infty}$ describing $L_{\infty}$-algebras.
    Remind that a $L_{\infty}$-algebra $V$ is a complex of vector spaces
     endowed with a coderivation $d_C$ of degree $(+1)$
      of the cofree cocommutative coassociative
    $\Z$-graded  coalgebra without counit cogenerated by $V[1]$:
      $$C:=\oplus_{n\ge 1} \bigl((V[1])^{\otimes n}\bigr)_{S_n}$$
      such that the component of $d_C$ mapping $V$ to $V$ coincides
       (up to a shift) with the differential $d_V$ of $V$.
Analogously, the operad $\mathsf{Assoc}$ is quasi-isomorphic to the operad
$\mathsf{Chains}(C_1)$, and also to the operad $A_{\infty}$ responsible for 
$A_{\infty}$-algebras.
   
\subsection{Formality of chain operads}    
 
\begin{theorem} 
The operad $\mathsf{Chains}(C_d)\otimes \R$ of complexes of real vector 
spaces is quasi-isomorphic to its cohomology operad endowed with zero
differential.
\end{theorem}
  
  In general, differential graded algebras quasi-isomorphic to their
cohomology endowed with zero differential, are called formal.
  The classical example is the de Rham complex of a compact K\"ahler
manifold. The result of Deligne-Griffiths-Morgan-Sullivan (see [DGMS])
is that  this algebra is formal as differential graded commutative 
associative algebra. Our theorem says that $\mathsf{Chains}(C_d)\otimes \R$ 
is formal as an algebra over the colored operad $\mathsf{Operads}$. 
   
   The story of this theorem is truly complicated. It was stated in the
preprint of Getzler and Jones originally for the cases $d=1$ and $d=2$. 
The case $d=1$ is trivial (as well as $d=0$). The authors referred to me 
for the claim that $\mathsf{Chains}(C_d)\otimes \R$ is not formal for
$d\ge 3$. Later it became clear that the proof in the Getzler-Jones preprint
is not correct. Personally, I thought for several years that even the fact 
is not true, and made several times calculations demonstrating the 
non-formality of $\mathsf{Chains}(C_2)\otimes \R$. The story began to
stir again when D.~Tamarkin posted on the net a new proof of the formality 
theorem in deformation quantization. The basic intermediate result 
in Tamarkin's approach was the formality of $\mathsf{Chains}(C_2)\otimes \R$ 
for which he claimed a new (with respect to [GJ]) proof.
Unfortunately, there were several mistakes in  Tamarkin's  proof at that time
as well, until a new corrected version was posted (see [T1]) which did
not use Deligne's conjecture directly. In [T1] is proven that $H_*(C_2)$ 
acts on the Hochschild complex of any algebra $A$.

Also Tamarkin found a proof that $\mathsf{Chains}(C_2)\otimes \C$
is quasi-isomorphic to $H_*(C_2)\otimes \C$ (see [T2]). Deligne's conjecture 
follow from the combination of two results of Tamarkin, but unfortunately 
not in a purely topological/combinatorial way. 
In the meantime I found a flaw in my calculations which ``show''
non-formality of $\mathsf{Chains}(C_2)\otimes \R$, and found a new proof 
of the formality of $\mathsf{Chains}(C_d)\otimes \R$ valid for all $d$.
The quasi-isomorphism which I constructed differs essentially from the 
one constructed by Tamarkin in [T2]. It seems that it gives  really 
different deformation quantizations of Poisson manifolds. The question of
differences between  quantizations is addressed in the third part of
this article.
 
I finish this subsection by a description of the homology operads of
topological operads $C_d$.
    
\begin{theorem} 
An algebra over $H_*(C_d)$ is 

\noindent  1) a complex in the case $d=0$,
  
\noindent  2) a differential graded associative algebra in the case $d=1$,
  
\noindent  3) a differential graded ``twisted'' Poisson algebra with 
the commutative associative  product of degree $0$ and 
  with the Lie bracket of degree $(1-d)$ in the case of odd $d\ge 3$,
   
\noindent   4) a differential graded ``twisted'' Gerstenhaber algebra
 with the commutative associative product of degree $0$ and with the 
Lie bracket of degree $(1-d)$ in the case of even $d\ge 2$.
   
   In cases 3) and 4) the Lie bracket satisfies the Leibniz rule with respect
    to the product. 
\end{theorem}
   
  The ``twisting'' in cases 3) and 4) means just that the commutator
   has the usual $\Z/2\Z$-grading, but a \textit{different} $\Z$-grading.
   The total rank of the homology group $H_*(C_d(n))$ is $n!$ for all $d\ge
1$. The rank of the top degree homology group of $C_d(n))$ is $(n-1)!$ 
for any $d\ge 2$ and $n\ge 1$. As a representation of $S_n$ this 
homology group  $H_{(n-1)(d-1)}(C_d(n))$ is isomorphic (up to tensoring 
by the sign representation of $S_n$ for even $n$) to $\mathsf{Lie}(n)$, 
the $n$-th space of the Lie operad.
     
The conclusion from two theorems in this section is that $d$-algebras for
$d\ge 2$ are essentially the same as twisted Poisson or Gerstenhaber 
algebras (depending on the  parity of $d$).    

\subsection{Sketch of the proof of the formality of chain
            operads}    

 The proof presented here is quite technical and it is not really essential
 for the rest of this paper. The reader can skip it and go directly to
section 3.4.

\subsubsection{Fulton-MacPherson compactification}    

Let us fix dimension $d\ge 1$. There is a modification $FM_d$ of 
the topological operad $C_d$ which is more convenient to work with. The
idea to consider operad $FM_d$ was proposed by several people, 
in particular by Getzler and Jones in [GJ]. The letters ``FM'' stand
for Fulton and MacPherson, who introduced  closely related constructions 
in the realm of algebraic geometry (see [FM]). 
I shall describe now the operad $FM_d$.

For $n\ge 2$ denote by $\tc_d(n)$  the quotient space of
  the configuration space of $n$ points in $\R^d$
  $$\mathsf{Conf}_n(\R^d):=\{(x_1,\dots,x_n)\in (\R^d)^n| 
\,x_i\ne x_j \,\,\,{\mathrm{for\,\,\,any}}\,\,\,i\ne j\}$$
   modulo the action of the group $G_d=\{x\mapsto \lambda x +v|\, \lambda \in 
  \R_{>0},\,\,v\in \R^d\}$. The space $\tc_d(n)$ is a smooth manifold 
of dimension $(nd -d-1)$. For $n=2$, the space $\tc_d(n)$ coincides with the
$(d-1)$-dimensional sphere $S^{d-1}$. There is an obvious free action of 
$S_n$ on $\tc_d(n)$. We define the spaces $\tc_d(0)$ and $\tc_d(1)$ to be 
empty. The collection of spaces $\tc_d(n)$ does not form an operad because 
there is no identity element, and compositions are not defined. 
    
The components of the operad $FM_d$ are
   
   1) $FM_d(0):=\emptyset$,
   
   2) $FM_d(1)=$point, 
   
   3) $FM_d(2)=\tc_d(2)=S^{d-1}$,
   
   4) for $n\ge 3$ the space $FM_d(n)$ is a manifold with corners, its
interior is $C'_d(n)$, and all boundary strata are certain products 
of copies of $\tc_d(n')$ for $n'<n$.
    
    A manifold with corners looks locally as a product of manifolds with
boundary. Any manifold with corners is automatically a \textit{topological}
manifold, although its \textit{smooth} structure is not the one of 
a differentiable manifold with boundary.
      
Intuitively, if a configuration of $n$ points in $\R^d$ moves in such a
way that several points (or groups of points) become at the limit close to 
each other, we use a microscope with a very large magnification
(apply a large element of the group $G_d$) in order to see in details 
the shape of configurations of points in clusters.
One of the rigorous definitions of $FM_d(n)$ is the following:
     
\begin{definition} 
For $n\ge 2$, the manifold with corners $FM_d(n)$ is the closure
of the image of $\tc_d(n)$ in the compact manifold 
$\left(S^{d-1}\right)^{n(n-1)/2}$ under the map
$$G_d\cdot (x_1,\dots,x_n)\mapsto
\left({x_j-x_i\over |x_j-x_i|}\right)_{1\le i<j\le n}$$
\end{definition}
      
Set-theoretically, the operad $FM_d$ is the same as the free operad 
generated by the collection of sets $(\tc_d(n))_{n\ge 0}$ endowed with  
$S_n$-actions as above.
     
It is possible to define another topological operad $FM'_d$, and 
two homomorphisms of operads
$$f_1:C_d\ra FM'_d,\,\,\, f_2:FM_d\ra FM'_d$$
such that for all $n\ge 0$ maps
$$f_1(n):C_d(n)\ra FM'_d(n),\,\,\, f_2(n):FM_d(n)\ra FM'_d(n)$$
are homotopy equivalences. In a sense, the spaces $FM'_d(n)$ parametrize
configurations of small disks in the standard disk, together
with a class of ``degenerate'' configurations in which some
(or all) disks are infinitely small.
We leave as an exercise to the reader to give a complete definition of $FM'_d$.
            
Applying the functor $\mathsf{Chains}$ we get two quasi-isomorphisms
of operads of complexes 
             $$\mathsf{Chains}(C_d)\ra \mathsf{Chains}(FM'_d)$$
             $$\mathsf{Chains}(FM'_d)\la \mathsf{Chains}(FM_d)$$

\subsubsection{A chain of quasi-isomorphisms}    
    
 For any $d\ge 2$  I shall define several operads of complexes 
 and construct a chain of quasi-isomorphisms between them. 
The total diagram is the following:
  $$\mathsf{Chains}(FM_d) \la \mathsf{SemiAlgChains}(FM_d)$$
\begin{equation}
\mathsf{SemiAlgChains}(FM_d)\otimes \R \ra \mathsf{Graphs}_d {\hat{\otimes}}\R
\end{equation}
$$\mathsf{Graphs}_d \la \mathsf{Forests}_d\la H_*(\mathsf{Forests)}H_*(C_d)$$

I shall now explain the first line. It is really technical, and is
introduced only to circumvent some difficulties with integrals which appear
in the second line (see the next section).
 
The operad $\mathsf{SemiAlgChains}(FM_d)$ is the suboperad of 
$\mathsf{Chains}(FM_d)$ consisting of combinations of maps
  $[0,1]^k\ra FM_d(n)$ whose graphs are real semi-algebraic sets. This is
well defined because the space  $FM_d(n)$ can be  described in terms of 
algebraic equations and inequalities. It is a closed semi-algebraic subset
in the product of $\frac{1}{2}n(n-1)$ copies of $S^{d-1}$. The natural 
inclusion of semi-algebraic chains into all continuous chains is
a quasi-isomorphism of operads.

\subsubsection{Admissible graphs and corresponding differential forms}    
\medskip

\begin{definition} 
An admissible graph with parameters $(n,m,k)$ (for $n\ge 1$ and $m\ge 0$)
 is a finite graph $\G$ such that 
                                                           
\noindent 1) it has no multiple edges,
 
\noindent  2) it contains no simple loops (edges connecting a vertex 
with itself),
 
\noindent  3) it contains $n+m$ vertices, numbered from $1$ to $n+m$,
 
\noindent  4) it contains $k$ edges, numbered from $1$ to $k$,
 
\noindent  5) any vertex can be connected by a path with a vertex whose 
index is in $\{1,\dots,n\}$,
 
\noindent  6) any vertex with index in $\{n+1,\dots,n+m\}$ has valency 
(i.e. degree) $\ge3$.
greater than or equal to $3$.

\noindent 7) for every edge $E$ of $\G$ we choose an orientation of this edge,
   i.e. we order the 2-element set of vertices to which $E$ is attached.

\noindent 8) if $n=1$ then the graph consists of just one vertex and has 
no edges, its parameters are $(1,0,0)$.
\end{definition}
   
Notice that although we endowed edges with orientations in 7),
we treat in 1)-6) the graph $\G$ as an \textit{unoriented} graph.
The structure of an admissible graph is completely determined by the 
attachment map 
$$\{1,\dots,k\}\ra \{(i,j)|1\le i,j,\le n+m,\,\,i\ne j\}$$
from the set of edges to the set of ordered pairs of distinct vertices.

\begin{definition}
  Let $\G$ be an admissible graph. We define $\omega_{\G}$ to be the 
differential form on $FM_d(n)$  given by the formula
\begin{equation}
\omega_{\G}:=(\pi_1)_*\circ \pi_2^*\left(\bigwedge_{{\mathrm{edges\,\,of\,\,}}
\G} \mathrm{Vol}_{S^{d-1}}\right)
\end{equation}
where
$$\pi_1:FM_d(n+m)\ra FM_d(n)$$ is the natural map, defined by forgetting 
the last $m$ points in the configuration of $(n+m)$ points in $\R^d$,
   $$\pi_2:FM_d(n+m)\ra \left( FM_d(2)\right)^k$$
   is the product of forgetting maps $FM_d(n+m)\ra FM_d(2)=S^{d-1}$
    associated with the edges of $\G$ (i.e. with ordered pairs of indices in
     $\{1,\dots,n+m\}$),
     $$Vol_{S^{d-1}}\in \Omega^{d-1}(S^{d-1})$$ denotes the volume form on
$S^{d-1}$ invariant under the action of the rotation group $SO(d,\R)$ 
and normalized such that the total volume 
$\int_{S^{d-1}}\mathrm{Vol}_{S^{d-1}}$ is $1$.
\end{definition}
  
  The degree of the form $\omega_{\G}$ is equal to
$$(d-1)k-dm=\mathrm{dim}\left(FM_d(2)^k\right)-\left(\mathrm{dim}(FM_d(n+m)-
\mathrm{dim}(FM_d(n))\right).$$ 
  If one changes orientations of edges, or the enumeration of edges,
    or the enumeration of vertices with indices from $n+1$ to $n+m$, one 
obtains the same form up to a sign. We leave as an easy exercise to the 
reader to write an explicit formula for this sign.
    
\begin{definition}
For every $n\ge 1$ we define $\mathsf{Graphs}_d(n)$ to be the 
$\Z$-graded vector space of all $\Q$-valued functions on 
the set of equivalence classes of admissible graphs with parameters
$(n,m,k)$, $m$ and $k$ arbitrary, such that if we change 
the enumeration of edges, or the enumeration of vertices with 
indices from $n+1$ to $n+m$, then the value of the function will be
 multiplied by  an appropriate sign as explained above.
 We define  $\Z$-grading of a function concentrated on one given 
equivalence class $[\G]$ with parameters $(n,m,k)$ as
 $\left(dm-(d-1)k\right)$.
\end{definition}
    
Unfortunately, the forms $\omega_\G$ are not $C^{\infty}$-forms on the
boundary of $FM_d(n)$. Still, for any $\G$ the integral of $\omega_{\G}$
over any semi-algebraic chain is absolutely convergent because the 
calculation of this integral reduces to the calculation the  volume form
over  a compact semi-algebraic chain of the top degree in the product 
of spheres. The total volume is finite because the multiplicity of a 
semi-algebraic map is bounded. Thus the integral of the volume form is
 convergent. 
         
{}From this follows that the form $\omega_\G$ gives a well-defined functional
on $\mathsf{SemiAlgChains}(FM_d)$. Turning things around, we can say
that any semi-algebraic chain gives a functional on the set
of equivalence classes of admissible graphs, i.e. an element
of the real-valued version of $\mathsf{Graphs}_d(n)$, the completed tensor
product $\R\hat\otimes \mathsf{Graphs}_d(n)$. The difference between the 
completed and the usual tensor product by $\R$ will eventually disappear 
because we shall meet the operad $H_*(FM_d)$ whose components 
are \textit{finite-dimensional}.

\begin{lemma} 
For any $\G$ the form $d\omega_\G$ (considered as a functional on 
semi-algebraic chains) is equal to the sum with appropriate signs of 
forms $\omega_{\G'}$ where the admissible graph $\G'$ is obtained
from $\G$ by contraction of one edge.
\end{lemma}

This lemma was proven in [K1]  (lemma in 6.6) for the case $d=2$ and 
in [K2] (lemma 2.1) for the case $d\ge 3$.
       
Thus, the graded space $\mathsf{Graphs}_d(n)$ of functions on graphs 
carries a naturally defined differential, and forms a complex of vector 
spaces over $\Q$. Moreover, the restrictions of the forms $\omega_\G$ to 
irreducible components of the boundary of the manifolds $FM_d(n)$ are 
finite linear combinations of products of analogous forms for simpler graphs.
This means that $\mathsf{Graphs}_d$ is an operad of complexes,
and the integration defines a homomorphism of operads of complexes of 
real vector spaces
$$\R {\otimes}\mathsf{SemiAlgChains}(FM_d)\ra
\R {\hat\otimes}\mathsf{Graphs}_d $$
           
It is not obvious that this arrow is a quasi-isomorphism.
This follows from an explicit calculation of the cohomology operad of the 
graph operad which we perform in the next subsection.
 
\subsubsection{Forests and Tree complexes}    
\medskip

\begin{definition} 
An admissible graph  is a forest iff it contains no nontrivial closed paths.
\end{definition}

 It is easy to see that for non-forest graph $\G$ the differential 
$d\omega_{\G}$ is a linear combination of forms $\omega_{\G'}$ for 
non-forest graphs $\G'$. Also, the restriction  of $\omega_{\G}$ to 
irreducible components of the boundary of  $FM_d(n)$ is a linear 
combination of products $\omega_{\G_{1}}\times\omega_{\G_{2}}$ 
where at least one of smaller graphs $\G_1,\G_2$ is not a forest.
All this implies that the following definition gives an operad:

\begin{definition}     
 The Operad $\mathsf{Forests}_d$ is a suboperad of $\mathsf{Graphs}_d$
  consisting of functions vanishing on all non-forest graphs.
\end{definition}
  
   A simple spectral sequence shows that the embedding
   $$\mathsf{Forests}_d\ra \mathsf{Graphs}_d$$ 
  is a quasi-isomorphism. We shall not write it here explicitly, only
mention the main idea of the calculation. Any graph which is not a forest
 contains a nonempty maximal subgraph $\G_{\mathrm{core}}$ such that 
the valency (in $\G_{\mathrm{core}}$) of each vertex of $\G_{\mathrm{core}}$ 
is at least $2$. The graph $\G$ is obtained from $\G_{\mathrm{core}}$ by 
attaching several trees to vertices of $\G_{\mathrm{core}}$, and also a forest
not connected with $\G_{\mathrm{core}}$. The desired result follows from
the vanishing of  the first term of the spectral sequence associated with
the filtration by  the number of vertices in $\G_{\mathrm{core}}$ on 
the quotient complex  ${\mathsf{Graphs}}_d(n)/{\mathrm{Forests}}_d(n)$.

Any admissible graph $\G$ with parameters  $(n,m,k)$  defines
a partition of $\{1,\dots,n\}$ into pieces corresponding to 
connected components of $\G$. All the graphs $\G'$ which appear in the 
decomposition of $d\omega_{\G}$ (as in the lemma in 3.3.3) give the same 
partition of $\{1,\dots,n\}$ as $\G$. This implies that the forest
complex ${\mathsf{Forests}}_d(n)$ splits naturally for any $n$ into 
a direct sum of subcomplexes corresponding to partitions 
 (i.e. equivalence relations) of the set $\{1,\dots,n\}$. 
     
The subcomplex of $\mathsf{Forests}_d(n)$ associated with a partition
is the tensor product of tree complexes over pieces of this partition, 
where \textit{trees} are non-empty connected forests, as usual.
    
We show an example of the tree complex in the case of $3$ vertices.
 There are $4$ admissible trees with the parameter $n$ equal to $3$
  (up to changing the enumeration of edges
  and orientations of edges):

\begin{center}
\scalebox{1}{\includegraphics{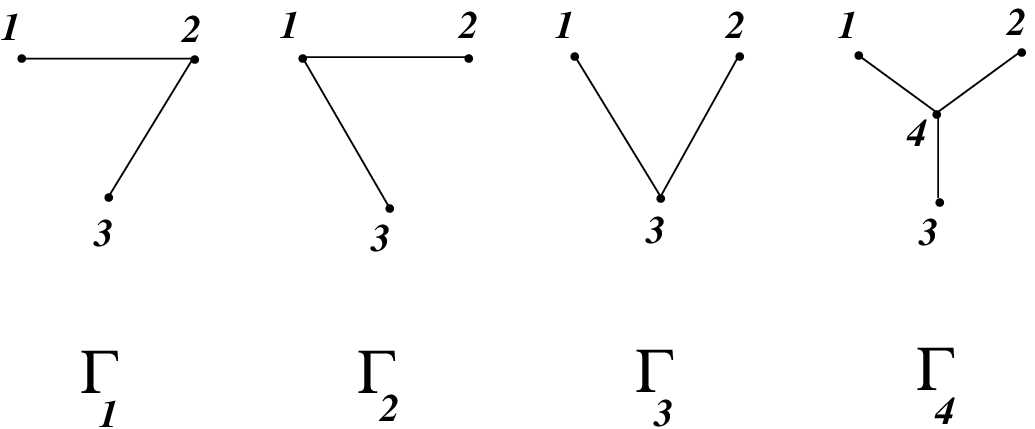}}
\end{center}

The formula for the de Rham differential of forms $\omega_{\G_i}$ is
$$d\omega_{\G_1}=d\omega_{\G_2}=d\omega_{\G_3}=0,\,\,\,d\omega_{\G_4}
=\pm\omega_{\G_1}\pm\omega_{\G_2}\pm\omega_{\G_3}\,\,\,.$$

The tree complex is the \textit{dual} complex to the complex spanned by
the differential forms $\omega_{\G}$ for trees $\G$. It is easy to see that 
the tree complex has the non-zero cohomology space only in the lowest 
degree. This implies that there is a natural map from the cohomology of the 
tree complex to the tree complex. Thus, we get a natural morphism from 
the cohomology of the forest complex to the forest complex. 
One can check that this map is a morphism of operads, and we get a
quasi-isomorphism $$H^*(\mathsf{Forests}_d)\ra \mathsf{Forests}_d$$
Comparing the cohomology of the forest complex with known results
on the cohomology of configuration spaces, we see that
        $$H^*(\mathsf{Forests}_d)=H_*(C_d)=H_*(FM_d)\,\,\,.$$
Because we proved already that 
$H^*(\mathsf{Graphs}_d)=H^*(\mathsf{Forests}_d)$,
we constructed a chain of quasi-isomorphisms as promised.

\subsection{Application to deformation quantization}    

The theorems in 3.2 show that any Gerstenhaber algebra, i.e.
an algebra over $H_*(C_d)$, can be canonically endowed with a $d$-algebra
structure. We prove in this section the result of Tamarkin:
 
\begin{theorem} 
Let $A:=\R[x_1,\dots,x_n]$ be the algebra of polynomials
  considered just as an associative algebra.
   Then the Hochschild complex $\mathsf{Hoch}(A)$ is quasi-isomorphic
  as $2$-algebra to its cohomology 
  $$B:=H^*(\mathsf{Hoch}(A))={\rm space\,\,of \,\,polynomial\,\,
   polyvector\,\,fields\,\,on\,\,\,}\R^n$$
   considered as a Gerstenhaber algebra, hence a $2$-algebra.
\end{theorem}

 As was mentioned in the introduction, this theorem implies
  the formality theorem from [K2]. 
  
  It is well-known (Hochschild-Kostant-Rosenberg theorem) that
   the space of polyvector fields is equal to the cohomology of 
$\mathsf{Hoch}(A)$, and the cup-product and the Lie bracket on 
$B=H^*(\mathsf{Hoch}(A))$ are the usual cup-product and the 
Schouten-Nijenhuis bracket on polyvector fields, respectively.
{}From the general formalism of  deformation theory it follows that
if $\mathsf{Hoch}(A)$ is \textit{not} quasi-isomorphic to $B$ then
there will be a \textit{non-zero} obstruction
element $\gamma_{\mathrm{obstr}}$  in the first cohomology group 
of the deformation complex $\mathsf{Def}(B)$ of $2$-algebra $B$. 
Moreover,  there exists an $\mathsf{Aff}(\R^n)$-invariant splitting of 
$\mathsf{Hoch}(A)$ into the direct sum of $B$ and a splitted contractible
complex (see arguments in [K1], 4.6.1.1). By the lemma in 2.6, this 
splitting induces a structure of $2$-algebra on $B$, and this structure 
is $\mathsf{Aff}(\R^n)$-invariant. It implies that 
$\gamma_{\mathrm{obstr}}\in H^1(\mathsf{Def}(B))$ is also
$\mathsf{Aff}(\R^n)$-invariant.
       
Let us calculate the cohomology of the deformation complex of the
$2$-algebra $B$. As a first approximation, we calculate the Hochschild 
cohomology of $B$:
  
\begin{theorem} 
The Hochschild complex of $B$ is quasi-isomorphic to $\R^1$ 
placed in degree $0$.
\end{theorem}
 
 In order to prove it we need an additional result concerning the
  Hochschild cohomology of $d$-algebras with trivial Lie bracket.
   By results from 3.2 we can speak about algebras over $H_*(C_d)$ instead
    of $d$-algebras.
      
\begin{lemma} 
Let $d\ge 2$ be an integer. Consider the algebra of polynomials in a finite
number of $\Z$-graded variables
     $$\R[x_1,\dots,x_N],\,\,\,deg(x_i)=d_i\in \Z$$
as an algebra over the operad $H_*(C_d)$, endowing it with zero differential 
and with the vanishing Lie bracket. Then the Hochschild cohomology of 
this algebra is, as $\Z$-graded vector space, the same as the algebra of 
polynomials in the doubled set of variables
$$\R[x_1,\dots,x_N,y_1,\dots,y_N],\,\,\,deg(x_i)=d_i,\,\,deg(y_i)=d-d_i$$
\end{lemma}
        
The statement of this lemma is similar to the classical 
Hochschild-Kostant-Rosenberg calculation of the Hochschild cohomology 
of the algebra of polynomials considered as an associative algebra  
(i.e. as $1$-algebra). We shall not give here the proof of this lemma.   
           
In general, if $\O(X)$ is the algebra of functions on a smooth
$\Z$-graded algebraic super-manifold $X$, then the Hochschild
cohomology of $\O(X)$ considered as a $d$-algebra, coincides with
the algebra of functions on the total space of the twisted by $[d]$
 cotangent bundle to $X$:
              $$H^*(\mathsf{Hoch}(\O(X)))=\O(T^*[d] X)\,\,\,.$$

Applying the above lemma to the following $H_*(C_2)$-algebra with 
\textit{vanishing} Lie bracket 
    $$B_0:=\R[x_1,\dots,x_n,\xi_1,\dots,\xi_n],\,\,\deg(x_i)=0,\,
   \deg(\xi_i)=+1$$  
we get the Hochschild cohomology 
    $$C_0:=\R[(x_i),(\xi_i),(\eta_i),(y_i)]_{1\le i\le n},$$
    $$\deg(x_i)=0,\, \deg(\xi_i)=\deg(\eta_i)=
+1,\,\deg(y_i)=+2\,\,\,.$$
 Our 2-algebra $B$ is obtained from $B_0$ by switching on the Lie bracket:
    $$[x_i,x_j]=[\xi_i,\xi_j]=0,\,\,\,[x_i,\xi_j]=\delta_{ij}\,\,.$$
It is easy to see that this produces the following differential on $C_0$:
     
     $$d_1:=\sum_{i=1}^n \bigl(\eta_i{\partial \over \partial x_i}
     +y_i{\partial \over \partial \xi_i}\bigr)$$
     The cohomology of the complex $(C_0,d_1)$ is equal to the de Rham 
cohomology of $\R^n$, i.e. to $\R$ placed in degree $0$. What we 
calculated is only the first term in a spectral sequence, but it is clear 
that higher differentials are zero because there is no space for them. 
This proves our theorem on the Hochschild cohomology of 
$B=\mathsf{Hoch}[x_1,\dots,x_n]$:
\begin{equation}
 H^*(\mathsf{Hoch}(\mathsf{Hoch}(\R[x_1,\dots,x_n])))=\R\,\,\,.
\end{equation}
     
     Now we are ready to calculate the first cohomology of the deformation
      complex of $B$. Remind (see 2.7) that there is a short exact sequence
      $$0\ra B[1]\ra \mathsf{Def}(B)\ra \mathsf{Hoch}(B)[2]\ra 0\,\,\,.$$
    Passing from this exact sequence to the level of the cohomology 
we get a long exact sequence 
$$\cdots\ra H^{i+1}(\mathsf{Hoch}(B))\ra B^{i+1}\ra 
H^i(\mathsf{Def}(B))\ra H^{i+2}(\mathsf{Hoch}(B))\ra \cdots$$
    which for $i=1$ says that 
    $$H^1(\mathsf{Def}(B))=B^2=\G(\R^d,\wedge^2 T)=
\{ {\rm \,polynomial \,\,bivector\,\,fields\,\,on\,\,}\R^d\}$$
    
   Then the argument goes as was sketched in the introduction.
   The explicit quasi-isomorphism between the space of polyvector fields
    and the Hochschild complex of the associative algebra
$A=\R[x_1,\dots,x_n]$ can be made invariant under the action of the group of 
affine transformations. There is no non-zero $\mathsf{Aff}(\R^n)$-invariant
bi-vector fields on $\R^n$. Thus the first non-trivial obstruction 
$\gamma_{\mathrm{obstr}}$ to the existence of a quasi-isomorphism between 
$B$ and $\mathsf{Hoch}(A)$ cannot exist, and $B$ is quasi-isomorphic to 
$\mathsf{Hoch}(A)$.
  
\section{Grothendieck-Teichm\"uller group action on
quantizations}    

\subsection{Integrals in deformation quantization}

There are several situations in deformation quantization where coefficients
in formulas are given by  explicit integrals.

 A) In Drinfeld's study of classical Knizhnik-Zamoldchi\-kov equations 
(see [D]), a formal  series in two non-commuting variables 
  appears, called an \textit{associator}. The coefficients in the explicit 
formula for the associator are  iterated integrals
\begin{equation}
I_{\eps_1,\dots,\eps_n}:=\int_{0<t_1<\dots<t_n<1} \omega_{\eps_1}(t_1)
 \wedge\dots\wedge\omega_{\eps_n}(t_n)
\end{equation}
 where $\eps_i\in\{0,1\},\,\,\,\eps_0=1,\eps_n=0,\,\,\,
\omega_0(t)=dt/t,\,\,\,\omega_1(t)=dt/(1-t)$.
 
 Among these numbers there are values of the Riemann zeta-function at positive
  integers, $\zeta(n)=I_{1,0,0,\dots,0}$ for a $n$-dimensional integral. 
In general, the integrals $I_{\eps_1,\dots,\eps_n}$ can be identified with
   so called \textit{multiple zeta-values} (see [Z]).
 
B) The same class of numbers appears in the Etingof-Kazhdan quantization of
Poisson-Lie algebras because it was based on Drinfeld's work (see [EK]). 
  
C) Tamarkin  uses the Drinfeld associator and Etingof-Kazhdan results 
from [EK] in his proof of the formality of $\mathsf{Chains}(C_2)\otimes \C$ 
(see [T2]).
 
 D) Tamarkin uses also the Drinfeld associator in the new proof of the 
formality theorem in deformation quantization (see [T1]). 
  
E) In my construction of quantization of Poisson manifolds (see [K1]) 
other types of integrals were used. These integrals are real-valued and 
are expressed by:

\begin{equation}
\int_{U_{n,m}} \bigwedge_{k=1}^{m+2n-2}\alpha(z_{i_k},z_{j_k}),\,\,m\ge 0, 
 n\ge 1
\end{equation}
 where the domain of integration $U_{n,m}$ is
 \begin{eqnarray*}
&& \{(z_1,\dots,z_{n+m})\in\C^{n+m}|\,
 z_1,\dots,z_m\in \R, z_1<\dots< z_m;\\
&&\quad \Im(z_{m+1}),\dots,\Im(z_{m+n})>0;z_i\ne z_j;
 z_{m+n}=\sqrt{-1}\}
\end{eqnarray*}
 
The form $\alpha(z,w)$ is 
\begin{equation}
{1\over 4 \pi i } d\log\left({(z-w)(z-{\overline w})
 \over ({\overline z}-w)({\overline z}-{\overline w})}\right)\,.
\end{equation}
The range of indices $i_k\ne j_k$ is $\{1,\dots, n+m\}$ for all
$k\in\{1,\dots,m+2n-2\}$.
  
F) Recently I realized that one can modify the above formulas,
replacing $\alpha$ by 
\begin{equation}
\alpha_{new}(z,w):= d\log\left({(z-w)\over ({\overline z}-w)}\right)
\end{equation}
 and dividing the integral by $(2\pi i)^n$.
 In this way one gets complex-valued integrals, and all identities
 proven  in [K2] remain true.
 
G) The last example is the use of integrals in the proof of formality of 
 $\mathsf{Chains}(C_2)\otimes \R$ presented in section 3.3.

I claim that all these integrals are closely related, and their use is 
probably unavoidable.

\subsection{Torsors}    
 
   In all the above situations A)-G),  we are constructing
    \textit{identifications} (up to a homotopy)
     between certain pairs of algebraic structures. For example, in 
deformation quantization we construct an isomorphism in the  homotopy
     category of Lie algebras between  the shifted by $[1]$
     Hochschild complex of $\R[x_1,\dots,x_n]$, and 
     the graded Lie algebra of polyvector fields on $\R^n$.
     
 Let $C$ be any category and $\E,\F$ be two isomorphic objects
 in this category. The set of isomorphisms  $$\mathsf{Iso}(\E,\F)$$
 is a non-empty set on which the group $\mathsf{Aut}(\E)$ acts simply 
transitively. The group $Aut(\F)$ acts also simply transitively on it, 
and the two actions commute. One can encode all these structures in a 
single map
\begin{equation}
    \mathsf{Iso}(\E,\F)^3\ra
        \mathsf{Iso}(\E,\F),\,\,\,(a,b,c)\mapsto a\circ b^{-1}\circ c
 \end{equation}      
     
\begin{definition} 
A torsor is a non-empty set $X$ endowed with a map
      $X\times X\times X\ra X$ satisfying the same identities as 
       maps $(a,b,c)\ra a b^{-1} c $ in groups. 
\end{definition}

A torsor is the same as a principal homogeneous space over a group.
One can give a more transparent definition of essentially the same
        structure:

\begin{definition}
A torsor is a category $C$ with only two objects $\mathsf{Ob}_1$ and 
$\mathsf{Ob}_2$, such that all the morphisms in  $C$ are invertible and the
objects $\mathsf{Ob}_1$ and $\mathsf{Ob}_2$ are equivalent. 
\end{definition}      
       If $X$ is a torsor (by Definition (18)) then one has \textit{two}
 groups acting simply transitively on $X$. Any element $x\in X$
 gives an identification between these two groups.

\subsubsection{Pro-algebraic torsors}    

In each of the cases listed in 4.1 in order to construct
 an isomorphism one is brought to solve an infinite system of quadratic 
equations with integer coefficients. For example, in deformation
 quantization we need to choose weights for all finite graphs in order to
   get  a quasi-isomorphism  between two homotopy Lie algebras. This
implies that the torsor of isomorphisms should be considered
    as an infinite-dimensional algebraic variety over $\Q$, not only
     as a set. We denote the torsor of isomorphisms considered
       as an algebraic variety by
         $${\underline {\mathsf{Iso}}}(\E,\F)$$
 Each of the two groups acting on the torsor is a projective limit 
of finite-dimensional affine algebraic groups over $\Q$.
The algebra of functions $\O(T)$ on a pro-algebraic torsor $T$
is a generalization of a Hopf algebra. Namely, it is a commutative 
associative unital algebra over $\Q$ together with the structure map 
\begin{equation}
 \O(T)\mapsto \O(T)\otimes \O(T)\otimes \O(T)
\end{equation}
 which is a homomorphism of algebras, and satisfies relations dual to 
the defining relations for set-theoretic torsors.
         We call this map the \textit{triple coproduct} in $\O(T)$.
         
The integrals in situations A)-G) from 4.1 give solutions to systems
 of quadratic equations in complex numbers, i.e. they give a homomorphism
 of algebras over $\Q$          $$\O(T)\ra \C\,\,\,$$
 where the pro-algebraic torsor $T$ depends on the concrete problem
 which we consider. In the next subsection we are going to describe a 
countable-dimensional commutative algebra $P$ over $\Q$ such that the 
homomorphism from $\O(T)$ to $\C$ is the composition of homomorphisms
       $\O(T)\ra P$ and $P\ra \C$.
      
\subsection{Periods, motivic Galois group, motives }    

Periods are integrals of algebraic differential forms with algebraic
 coefficients. The following numbers are periods:
 $$ 1,\,\,\sqrt{2},\  i=\sqrt{-1},\ \pi,\  \log(2),$$
 $$\zeta(3)=\int_{0<t_1<t_2<t_3<1}
 {dt_1\over 1-t_1} {dt_2\over t_2}{dt_3\over t_3},\
 \mathrm{elliptic\,\, integral}\ \int_1^2\sqrt{x^3+1}dx,\dots$$
 
 We shall need a more precise definition of periods. 
Let $X$ be a smooth
 algebraic variety of dimension $d$
 defined over $\Q$, $D\subset X$ be a divisor  with normal crossings
  (i.e. locally $D$ looks like a collection of coordinate hypersurfaces),
   $\omega \in \Omega^d(X)$ be an  algebraic differential
   form  on $X$  of top degree ($\omega$ is automatically closed), 
   and $\gamma\in H_d(X(\C),D(\C);\Q)$ be a (homology class of a)
    singular chain on the complex manifold $X(\C)$ with boundary on 
  the divisor $D(\C)$. With these data one associates the integral 
 $\int_\gamma \omega \in \C$.  We say that this number
  is the period of the quadruple $(X,D,\omega,\gamma)$. 
 One can always reduce convergent integrals of algebraic forms over 
 semi-algebraic sets defined over the field algebraic numbers $\overline{\Q}$
  to the form as above, using the functor of the restriction of scalars
   to $\Q$ and the resolution of singularities in characteristic zero.
    
    Usual tools for proving identities between integrals are
    the change of variables and the Stokes formula.
    Let us formalize them for the case of periods.

\begin{definition}
The space $P_+$ of effective periods is defined as a vector space over $\Q$ 
generated by the  symbols $[(X,D,\omega,\gamma)]$ representing equivalence
classes of quadruples as above, modulo the following relations:

 1) (linearity)  $[(X,D,\omega,\gamma)]$ is linear in both $\omega$ and
 $\gamma$
 
 2) (change of variables) If $f:(X_1,D_1)\ra (X_2,D_2)$
  is a morphism of pairs defined over $\Q$, 
  $\gamma_1\in H_d(X_1(\C),D_1(\C);\Q)$ and $\omega_2\in \Omega^d(X_2)$ then
  $$[(X_1,D_1,f^*\omega_2, \gamma_1)]=[(X_2,D_2,\omega_2,f_*(\gamma_1))]$$
  
  3) (Stokes formula) Denote by $\tilde{D}$ the normalization of $D$
   (i.e. locally it is the disjoint union of irreducible
   components of $D$), the variety  $\tilde{D}$ containing
 a divisor with normal crossing $\tilde{D}_1$
     coming from double points in $D$. If $\beta \in \Omega^{d-1}(X)$
      and $\gamma \in H_d(X(\C),D(\C);\Q)$ then
      $$[(X,D,d\beta,\gamma)]=[(\tilde{D},\tilde{D}_1, 
\beta_{|\tilde{D}},\partial\gamma)]$$
      where $\partial:H_d(X(\C),D(\C);\Q)\ra H_{d-1}(
      \tilde{D}(\C),\tilde{D}_1(\C);\Q)$ is the boundary operator.
\end{definition}
      
    It is conjectured in number theory that the evaluation homomorphism
    $P_+\ra \C$ is a monomorphism, i.e. all identities between periods
      can be proved using standard rules only. For example, the fact
      that number $\pi$ is transcendental follows from this conjecture.
      
The effective periods form an algebra because the product of integrals
   is again an integral (Fubini formula). The field of algebraic numbers 
${\overline\Q}\subset \C$  can be considered as a subalgebra over $\Q$ of  
the algebra $P_+$. An algebraic number $x\in {\overline \Q}$ which solves  
a polynomial equation $P=0,\,\,\,P\in \Q[t]$, is the period of 
$0$-dimensional variety $X\subset \A^1_{\Q}$ defined by the equation $P=0$.
The number $x$ gives a complex point of $X$, i.e. a $0$-chain.
The standard coordinate $t$ on the affine line $\A^1_{\Q}$ gives a
zero-form after restriction to $X$ whose pairing with 
$x\subset X(\C)\subset \A^1_{\Q}(\C)=\C$ is tautologically  equal to $x$.
        
   It is convenient to extend the algebra of effective periods to a
   larger algebra $P$ by inverting formally the element whose evaluation in 
$\C$ is $2\pi i$. Informally, we can write that the whole algebra of periods
      $P$ is $P_+[(2\pi i)^{-1}]$.
      
 The algebra $P$ is an infinitely generated algebra over $\Q$, but as any
algebra it is an inductive limit of finitely-generated subalgebras.
 This means that $\mathsf{ Spec}(P)$ is a projective limit of 
finite-dimensional affine schemes over $\Q$. We claim that 
$\mathsf{Spec}(P)$ carries a natural structure of a pro-algebraic 
torsor over $\Q$.
     
The formula for the structure map  $\Delta:P\ra P\otimes P\otimes P$ 
can easily be written in terms of period matrices of individual algebraic 
varieties. Namely, let $(P_{ij})$ be the period matrix of an algebraic variety
consisting of pairings between classes running through
a basis in $H_*(X(\C),\Q)$ and a basis in $H^*_{\mathrm{de \,Rham}}(X)$.
More generally, one should consider homology and cohomology of pairs
of algebraic varieties over $\Q$. It follows from several results in algebraic
geometry that the period matrix is a square matrix with entries 
in $P_+$, and  determinant in $\sqrt{\Q^{\times}}\cdot (2\pi i)^{\Z_{\ge 0}}$.
This implies that the inverse matrix has coefficients in the extended
 algebra $P=P_+[(2 \pi i)^{-1}]$.
          
\begin{definition} 
The triple coproduct in $P$ is defined by
\begin{equation}
\Delta (P_{ij}):=\sum_{k,l} P_{ik}\otimes (P^{-1})_{kl}\otimes P_{lj}
\end{equation}
for any period matrix $(P_{ij})$.
\end{definition}
          
We show how to calculate the triple coproduct in a simple example.
  Let $X:=\A^1_{\Q}\setminus \{0\}$ be the affine line with deleted
  point $0$, and $D:=\{1,2\}\subset X$ be a divisor in $X$.
  The first homology group of pair (relative homology)  

 $$H_1(X(\C),D(\C);\Q)=H_1(\C\setminus\{0\},\{1,2\};\Q)$$
is two-dimensional and is generated by two chains: a small closed 
anti-clockwise oriented path around $0$, and the  interval $[1,2]$.
 The algebraic de Rham cohomology group $H^1_{\mathrm{de\,Rham}}(X,D)$
is also two-dimensional, and is generated by the $1$-forms $dt$ and $dt/t$ 
where $t$ is the standard coordinate on $X\subset \A_{\Q}^1$.
 The period matrix is 
\begin{equation}
\pmatrix{1 & \log(2)\cr 0 & 2 \pi i\cr}.
\end{equation}
                 
 {}From this one can deduce that
 $$\Delta(2\pi i)=2\pi i \otimes {1\over 2\pi i}\otimes 2\pi i \,\,,$$
 $$\Delta(\log(2))=\bigl(\log(2)\otimes  {1\over 2\pi i}\otimes 2\pi i\bigr)
                 -\bigl(1\otimes {\log(2)\over 2\pi i}\otimes 2\pi i
                 \bigr)+ \bigl(1\otimes 1\otimes \log(2)\bigr )\,\,\,.$$
 It is not clear why the definition given above is consistent, because
  it is not obvious why the triple coproduct preserves the defining relations
 in $P$. This follows  more or less automatically from the following result:
    
\begin{theorem}[M.~Nori]
       Algebra $P$ over $\Q$ is the algebra of functions on the pro-algebraic
        torsor  of isomorphisms between two 
       cohomology theories, the usual topological cohomology theory
       $$H_{\mathrm{Betti}}^*:X\mapsto H^*(X(\C),\Q)$$
        and the algebraic de Rham cohomology theory
         $$H_{\mathrm{de\,Rham}}^*:X\mapsto {\bf H} ^*(X,\Omega^*_X)$$
\end{theorem}

The motivic Galois group in Betti realization $G_{M,\mathrm{Betti}}$ 
is defined as the pro-algebraic group acting on $\mathsf{Spec}(P)$ from the 
side of Betti cohomology. Analogously one defines the de Rham version 
$G_{M,\mathrm{de\,Rham}}$.  The category of motives is defined as the 
category of representations of the motivic Galois group. It does not matter 
which realization we choose  because the categories for both
    realizations can be canonically identified with each other.
    Here is the ``symmetric'' definition of the category of motives:
    
\begin{definition} 
The symmetric monoidal category of motives over $\Q$ with coefficients in $\Q$
is defined as the category of vector bundles on $\mathsf{Spec}(P)$ endowed
with two commuting actions of the motivic Galois groups 
$G_{M,\mathrm{Betti}}$ and $G_{M,\mathrm{de\,Rham}}$.
\end{definition}

     The following elementary definition gives a category equivalent 
      to the category of motives:
\begin{definition}
A framed motive of rank $r\ge 0$ is an invertible  $(r\times r)$-matrix
$(P_{ij})_{1\le i,j,\le r}$ with coefficients in the algebra $P$,
satisfying the equation
\begin{equation}
\Delta (P_{ij})=\sum_{k,l} P_{ik}\otimes (P^{-1})_{kl}\otimes P_{lj}
\end{equation}
for any $i,j$. The space of morphisms from one framed motive to another,
corresponding to matrices 
$$P^{(1)}\in GL(r_1,P),\,\,\,P^{(2)}\in GL(r_2, P),$$
        is defined as
        $$\{ T\in \mathsf{Mat}(r_2\times r_1,\Q)| \,T P^{(1)}=P^{(2)}T\}$$
\end{definition}
 
 The cohomology groups of varieties over $\Q$ can be considered as objects
  of the category of motives.  From comparison isomorphisms in algebraic 
geometry follows that there   are also $l$-adic realizations motives, 
on which the absolute   Galois group
    $\mathsf{Gal}({\overline{\Q}}/\Q)=\mathsf{Aut}({\overline{\Q}})$ acts.
    
     The usual Galois group $\mathsf{Gal}({\overline{\Q}}/\Q)$ is a
profinite group, projective limit of finite groups. It is the image of the 
motivic Galois group $G_{M,\mathrm{Betti}}$ in motives coming from $0$-th 
cohomology groups of schemes defined over $\Q$. It can be described also
 as follows: the subalgebra  ${\overline{\Q}}$ of $P$ is closed under
 the triple coproduct,
          $$\Delta_{|{\overline{\Q}}}
          :{\overline{\Q}}\ra {\overline{\Q}}\otimes_{\Q}{\overline{\Q}}
          \otimes_{\Q}{\overline{\Q}}$$
          and its spectrum is a quotient torsor of $P$ on which
           $\mathsf{Gal}({\overline{\Q}}/\Q)$ acts simply transitively.
          
   The definition of the category of motives given above
    is a natural generalization of the ``folklore'' 
    definition of the category of  pure motives given
     for example in [S]. There is an elaborated hypothetical picture
      of motives in number theory, from which follows that
      the category of  motives as defined here is exactly what is expected.
       The definition given here is equivalent
      to the one advocated by M.~Nori in [N]. Specialists in motives consider
       this definition as ``cheap'', and expect in the future
        something more elaborate and not directly referring to   
       explicit realization functors (like Betti realization etc.).
       
       The advantage of the definition given here is that it does not assume
    the validity of any conjecture, and is directly applicable to 
         the present study of operads in deformation quantization.
        
\subsection{Grothendieck-Teichm\"uller group }    
    
   All the integrals appearing in situations listed in 4.1 share several 
common features. The corresponding motives are so called 
\textit{mixed Tate motives}, which means in particular that the 
period matrix is upper-triangular in certain bases, with integral powers 
of $(2\pi i)$ on the diagonal (like in the example in the previous 
subsection). Also these motives are \textit{unramified over} 
$\mathsf{Spec}(\Z)$. This property can be expressed in terms of $l$-adic 
representations of $\mathsf{Gal}({\overline{\Q}}/\Q)$ and can be verified 
in concrete situations by checking that certain discriminants are equal
to $0$,$1$ or $-1$.  The period matrix (17) above, 
is \textit{not} unramified. Namely, it is ramified at prime $2$.
       
We denote by $P_{\Z,\mathrm{Tate}}$ the subalgebra of $P$ generated
by $(2\pi i)^{\pm 1}$ and by periods of mixed Tate motives unramified
 over $\mathsf{Spec}(\Z)$.  There is a conjectural picture for it:
               
\begin{conjecture}
The  quotient of the motivic Galois group (in the de Rham realization)
acting simply transitively on $\mathsf{Spec}(P_{\Z,\mathrm{Tate}})$  is 
 a pro-solvable connected group over $\Q$, an extension of the multiplicative
 groups scheme ${\bf G_m}=GL(1)$ by a pro-nilpotent
 group whose Lie algebra is free and generated by elements
 of weights $3,5,7,\dots$ (one element in each odd weight $\ge 3$)
 with respect to the adjoint action of ${\bf G_m}$.
\end{conjecture}

This conjecture follows from general Beilinson conjectures
on  motives and $K$-theory (see [Ne]).
           
There are two other conjectures concerning $P_{\Z,\mathrm{Tate}}$.
          
\begin{conjecture}
$P_{\Z,\mathrm{Tate}}$ is the subalgebra of $P$ generated by 
$(2\pi i)^{\pm 1}$ and by periods whose evaluation in $\C$ are integrals 
$I_{\epsilon_1,\dots, \epsilon_n}$ which appear in Drinfeld's associator 
(see 4.1.A).
\end{conjecture}             
       
  There are two reasons for this conjecture. First of all, if we believe 
in the picture of $P_{\Z,\mathrm{Tate}}$ explained above,
there are  $O(c^n)$ linearly independent effective unramified Tate motives
of weight $n$, where $c=1.32471...$ is the positive root of the equation
$x^3=x+1$. On the other hand, there are $2^{n-2}$ integrals 
$I_{\epsilon_1,\dots, \epsilon_n}$, and with high probability
(because $c<2$) there are enough integrals to  span the whole algebra 
$P_{\Z, \mathrm{Tate}}$. The second reason is also probabilistic,
computer experiments confirm that the Poincar\'e series of the algebra
 generated by integrals graded by weights 
$w(I_{\epsilon_1,\dots, \epsilon_n}):=n$ equals
to the expected series $1/(1-t^2-t^3)$ up to $O(t^{13})$.
          
The next conjecture concerns the so called Grothendieck-Teichm\"uller
group $GT$ (see [D]). 
           
\begin{conjecture} 
The quotient of the motivic Galois group
acting simply transitively on the spectrum of the subalgebra
of $P$ generated by $(2\pi i)^{\pm 1} $ and integrals
$I_{\epsilon_1,\dots, \epsilon_n}$, coincides with the group $GT$.
\end{conjecture}   
           
The group $GT$ is a pro-algebraic group over $\Q$, an extension of 
${\bf G_m}$ by a pro-nilpotent group. One of the definitions of $GT$ is 
as the group of automorphisms of the tower of pro-nilpotent completions
of pure braid groups $\pi_1(\C^n\setminus \mathrm{Diag})$ for all $n$.
The group of automorphisms of these pro-nilpotent completions
coincides with the group of automorphisms of the tower
of \textit{rational homotopy types} of classifying spaces of these groups. 
Classifying spaces are $(\C^n\setminus \mathrm{Diag})_{n\ge 2}$, 
configuration spaces of $\R^2=\C$. By general reasons the motivic Galois
group acts on the tower of rational homotopy types of 
$(\C^n\setminus \mathrm{Diag})=
\left(\A_{\Q}^n\setminus \mathrm{Diag}\right)(\C)$,
thus it maps to $GT$. Moreover, it is easy to see that the
periods which appear in the image of this action are exactly the
numbers $I_{\epsilon_1,\dots, \epsilon_n}$.
         
 A priori, there is no reason for the tower of Malcev completed
 braid groups to have a non-trivial automorphism.
Number theory provides (via motivic Galois group)
 a supply of such automorphisms, and the conjecture 4 above says
that there is no other automorphism.
             
The group $GT$ maps to the group of automorphisms  in homotopy sense
 of the operad $\mathsf{Chains}(C_2)$. Moreover, it seems to coincide with 
$\mathsf{Aut}(\mathsf{Chains}(C_2))$ when this operad is considered as an 
operad not of complexes but of differential graded cocommutative coassociative
coalgebras (strictly speaking there is no coproduct on singular
 cubical chains, but one can overcome this technical problem).
                 
 {}From now on we shall assume the validity of conjectures 2,3, and 4 
and identify the group acting simply transitively on $P_{\Z,\mathrm{Tate}}$
 with the Grothendieck-Teichm\"uller group $GT$.
             
\subsection{Conjectures about torsors}    
 
Remind that the integrals described in situations A)-G) in 4.1 give 
complex points of the corresponding torsors, i.e.  homomorphisms of algebras 
over $\Q$
 $$\O(T)\ra \C\,\,\,.$$
 Moreover, the values of the integrals in all these cases are periods, 
and we proved that quadratic equations are satisfied using only the 
Stokes formula. This implies that we have in fact a homomorphism
\begin{equation}
\O(T)\ra P_{\Z,\mathrm{Tate}}\,\,(\hookrightarrow P)
\end{equation}
   In terms of pro-algebraic affine schemes, this means that we have a map 
\begin{equation}
\mathsf{Spec}(P_{\Z,\mathrm{Tate}})=
{\underline {\mathsf{Iso}}}_{\Z,\mathrm{Tate}}
    (H_{\mathrm{de\,Rham}}^*,H_{\mathrm{Betti}}^*)\ra T
\end{equation}
    between two torsors considered just as pro-algebraic schemes.
     For example, in the formality theorem (cases D),E),F)) the torsor $T$
  is the torsor of isomorphisms between two homotopy Lie algebras
   $$T={\underline {\mathsf{Iso}}}(T^*_{\mathrm{poly}}, D^*_{\mathrm{poly}})$$
   (see [K1] for a precise definition of $T^*_{\mathrm{poly}}$ and 
$D^*_{\mathrm{poly}}$; one can safely replace $D^*_{\mathrm{poly}}$ by the 
shifted by $[1]$ Hochschild complex of the algebra
     $A:=\R[x_1,\dots,x_n]$, and $T^*_{\mathrm{poly}}$ by the algebra of 
polynomial polyvector fields on $\R^d$). It is natural to ask whether 
    the constructed map is in fact  a map of torsors.
        
\begin{conjecture} 
In cases A),B),C),D),F) the map from $P$ to the corresponding torsor of 
isomorphisms is a map of torsors.
\end{conjecture}

 I checked that the generator of $\Z/2\Z\subset G_{M,\mathrm{Betti}}$
 (acting via complex conjugation on $X(\C)$ where $X$
  is defined over $\Q$) corresponds to the natural involution
  of the Hochschild complex $C^*(A,A)$ for $A=\R[x_1,\dots,x_n]$,
  the  Hochschild complex being considered as a dg-Lie algebra. 
  The involution acts on $\phi\in C^n(A,A)$ as
             $$\phi\mapsto {\overline \phi},\,\,\,
             {\overline \phi}(a_1\otimes\dots \otimes a_n):=
             (-1)^{(n+1)(n+2)/2} \phi(a_n\otimes\dots \otimes a_1)\,\,\,.$$
             In terms of deformation theory this corresponds to the passage
              from a product $\star$ to the \textit{opposite} product 
               $$a{\overline \star} b:=b\star a$$
  Also, the group ${\bf {G_m}}$ acts on $\mathsf{Spec}(P_{\Z,\mathrm{Tate}})$ 
from the de Rham side, and also on 
${\underline {\mathsf{Iso}}}(T^*_{\mathrm{poly}}, D^*_{\mathrm{poly}})$
from the side of $T^*_{\mathrm{poly}}$ rescaling polyvector fields
according to $\Z$-grading. Again, I checked that the two actions of 
${\bf {G_m}}$ are compatible.
              
In cases the E),G) not listed in the conjecture, the values of the
integrals are \textit{real} numbers. It seems that the map between torsors
does not respect torsor structure. Moreover, I think that it is equal to 
the composition of a map of torsors with a certain universal
 map of schemes $P_{\Z,\mathsf{Tate}}\ra P_{\Z,\mathsf{Tate}}$
 which is \textit{not} a map of torsors.
               
\subsection{Incarnations of the  Grothendieck-Teichm\"uller group }    

I present here two examples in which one can see the action of $GT$ on
 deformation quantizations. Proofs will appear elsewhere.
 
\begin{theorem}
Let $\g$ be a finite-dimensional Lie algebra over $\R$.
     In the Duflo-Kirillov isomorphism (see [K1]) between the center
      of the universal enveloping algebra $U\g$, and the algebra of
      invariant polynomials on $\g^*$, one can replace the formal series 
      $$F(x)= \sqrt{e^{x/2}-e^{-x/2}\over x}$$
       by the product
       $$F_{\mathrm{new}}(x)=F(x)\cdot 
       \exp\left(\sum_{k=0}^{\infty} a_{2k+1} x^{2k+1}\right)$$
       where $a_1,a_3,a_5,\dots$ are arbitrary constants. 
       For any such choice one gets again an isomorphism 
       compatible with products.
\end{theorem}

       In particular, one can replace $F$ by the following entire function
        $$F_{\mathrm{nice}}(x)={1\over\G\bigl({x\over 2\pi i}+1\bigr)}\,\,\,$$
         for the choice 
         $$a_1={{\rm Euler \,\,\,constant}\over 2 \pi i}
         ={0.57721...\over 2\pi i},\,\,a_3={\zeta(3)
         \over 3(2\pi i)^3},\,a_5={\zeta(5)\over 5 (2\pi i)^5},\dots
         $$
        
        The Euler constant is (probably) not a period, and enters the
         formula only for \ae sthetical reasons. The terms with coefficients
          $a_3,a_5,\dots$ appear because of the action of $GT$.
          
        As a corollary to Theorem 7 above we have the following
        
\begin{theorem}
If $\g$ is a finite-dimensional Lie algebra,
then the differential operators on $\g^*$ with constant coefficients
which are Fourier transforms of the following polynomials on $\g$:
$$P_{2k+1}(\gamma):=\mathrm{Trace} (\mathrm{ad}(\gamma)^{2k+1}),\,\,\,k\ge 0$$
act on the subalgebra of $\mathrm{ad}(\g)$-invariant polynomials as
derivations, i.e. satisfy the Leibniz rule.
\end{theorem}

 Unfortunately, as I learned from M.~Duflo, for all finite-dimensional 
Lie algebras, operators such as in those considered in Theorem 7 above,
when restricted to $\mathsf{Sym}(\g)^{\g}$ are all equal to zero.
 Thus, we do not get anything visible here. Nevertheless, my result works 
also for  Lie algebras in  rigid symmetric monoidal categories, for example
for a finite-dimensional Lie superalgebra. There is a chance that one 
can find non-trivial examples there.
         
Another incarnation of the motivic Galois group is in a sense
Fourier dual to the previous one. Namely, let $X$ be a complex manifold 
(or a smooth algebraic variety over a field of characteristic zero).
       Define the Hochschild cohomology of $X$ as the following graded
       commutative associative algebra:
       $$HH^k(X):=\oplus_{i+j=k} H^i(X,\wedge^j(T_X))\,\,\,.$$
       The product on $HH^*(X)$ is given by the usual cup-product
        of polyvector fields and of cohomology classes.
        Every element of the \textit{Hodge} cohomology
        $$\oplus_{i,j} H^i(X,\wedge^j T^*_X)$$
         gives a linear operator on $HH^*(X)$. It comes from the 
         cup-product in cohomology and from the convolution operator
         $$\wedge^a T^*_X\otimes\wedge^b T_X\ra \wedge^{b-a} T_X$$
         acting fiberwise on the level of bundles.
         Operators corresponding to elements in $H^i(X,\wedge^i T^*_X)$
          act on $HH^*(X)$ preserving the $\Z$-grading.
          
      One can construct a version of characteristic classes
      of vector bundles on $X$ with values in the diagonal part
      $\oplus_{i\ge 0} H^i(X,\wedge^i T^*_X)$ of the Hodge cohomology. 
    This can be done using the Atiyah class of a vector bundle $\mathcal{E}$
 taking values in $H^1(X,T^*_X\otimes \mathsf{End}(\mathcal{E}))$. 
    This class is the class of the extension of $\mathcal{E}$ by 
   $T^*_X\otimes {\mathcal{E}}$ given by the bundle of $1$-jets of sections 
   of $\mathcal{E}$. The traces of powers of the Atiyah class give 
  characteristic classes associated with monomial symmetric functions, 
  and other characteristic classes are polynomials in basic ones.
                
\begin{theorem} 
Operators on $HH^*(X)$ corresponding to odd components 
        $$ch_{2k+1}(T_X)\in H^{2k+1}(X,\wedge^{2k+1}T^*_X)$$
    of the Chern character $ch(T_X)$ of the tangent bundle,
    are derivations of $HH^*(X)$ with respect to the cup-product.
\end{theorem}
                   
   The statement of this theorem is very easy for $k=0$, but
   I do not know   elementary proofs for higher values of $k$.
               
   In the case of $X$ being a compact Calabi-Yau variety
    one can identify the vector spaces $HH^*(X)$ and the usual 
cohomology $H^*(X)$ using a Hodge decomposition and the convolution
 with a volume element on $X$.  The operator corresponding 
 to the class $\mathrm{ch}_1(T_X)=c_1(T_X)=0$ is zero, but higher operators
    are in general non-zero. They correspond after the identification
     of  $HH^*(X)$ with   $H^*(X)$ to the multiplication operators
      by $\mathrm{ch}_{2k+1}(T_X)$.
      
      In general, the result is that the algebra  $HH^*(X)$ has
       an infinite set of commuting automorphisms  labeled by odd positive 
     integers, and it is a module over a solvable quotient of the
         motivic Galois group. The Lie algebra of this quotient has the
following  basis and bracket:
         $$L_0,Z_3,Z_5,Z_7,\dots,\,\,\,[L_0,Z_{2k+1}]=(2k+1)Z_{2k+1},\,\,\,
         [Z_{2k+1},Z_{2l+1}]=0\,\,\,.$$

\section{Relations to quantum field theories }    

\subsection{Local fields and $d$-algebras }    

There is no satisfactory definition yet of a quantum field theory (QFT).
 One expects at least that a QFT on a manifold $X$
  (``the space-time'') in the Euclidean framework
  gives a super vector bundle of ``local fields'' $\Phi$ on $X$
   (the bundle $\Phi$ may be infinite-dimensional), and  correlators 
  $$\langle\phi_1(x_1)\dots \phi_n(x_n)\rangle \in \C,\,\,\,n\ge 0$$
  which are even polylinear maps from the tensor product of fibers of
$\Phi$ at   pairwise distinct points
  $$\langle\dots\rangle: \Phi_{x_1}\otimes\dots\otimes \Phi_{x_n}\ra \C$$
  depending smoothly and $S_n$-equivariantly on
  $$(x_1,\dots,x_n)\in \mathsf{Conf}_n(X)=X^n\setminus \mathrm{Diag}$$
 Also one expects that there is an \textit{operator product expansion}
  of fields at points converging to one point.
  I am not going to discuss here in details the properties of the 
operator products expansion.
   
   It seems very plausible that with an appropriate definition 
    one can prove the following 

\begin{conjecture}
     For any QFT on vector space $V=\R^d$ with invariance under the action
of the group $G_d$ of parallel translations and dilatations
       (in particular, for any conformal field theory), on the tensor product
       $$\Phi_0\otimes \wedge^*(V^*)$$
       there is a structure of a $d$-algebra.
\end{conjecture}
 
       This implies that translation-invariant differential forms on   $\R^d$
       with values in local fields form 
        a  homotopy Lie algebra $\g^*$ whose graded components are
$$\g^k=\Phi_0\otimes \wedge^{k+d-1}(V^*),\,\,\,k\in \{-(d-1),\dots,+1\}.$$

        The moduli space associated with this homotopy 
         Lie algebra is unobstructed because $\g^2=0$,
    and it can be interpreted as the formal deformation theory of our QFT
    as a translation-invariant theory,
    with all renormalizations and regularizations automatically
    included in the structure of $L_{\infty}$-algebra on $\g^*$.
    Presumably, the construction of (higher) brackets
    on $\g^*$ is closely related with the Hopf algebra
    studied by Connes and Kreimer (see [CK]), where the
    case of free massless theory   is considered.
           
   Also, it seems that the notion of an action of $(d+1)$-algebra
   on a $d$-algebra is closely related with field theories
   on manifolds with boundaries. This subject became very popular
   in modern string theory after the Maldacena conjecture  
   on  boundary conformal field theories for QFT on anti-de-Sitter
   spaces (Lobachevsky spaces).
   Here I should notice that a work of Mosh\'e and co-authors
   (see [FFS] and references therein) 
   was one of the predecessors of the modern AdS picture.

\subsection{Action of the motivic Galois group on the moduli space of QFTs}  
         
  In the construction of perturbations
   of a given conformal field   theory one needs to  calculate
    (and often regularize) integrals of correlators over configuration spaces.
    In the free theory case the integrals are exactly Feynman integrals
     in the diagrammatic expansion.
     
     D.~Broadhurst and D.~Kreimer (see [BK])
     observed that  all Feynman diagrams up to 7 loops 
      in any QFT in \textit{even dimensions}  
      gives same numbers as appear in Drinfeld associator. It is not clear
a priori why this
       happens. In any case one can see immediately from formulas that
        all constants are in fact periods.
        
\begin{conjecture} The motivic Galois group $G_{M,Betti}$ 
       acts (in homotopy sense)  on the homotopy
       Lie algebra $\g^*$ associated with the free massless theory
        in any  dimension. In the case of even dimension the action
         factors through the quotient group $GT$ as in 3.4.
          The action should be somehow related with the action
         on values of Feynman integrals.
\end{conjecture}

         I have only one ``confirmation'' of this conjecture,
          that is the deformation quantization story. 
          In my proof of the formality theorem, the integrals which 
           appear in the explicit formula come from Feynman diagrams
         for a perturbation  of a free two-dimensional quantum field
         theory on the Lobachevsky plane (see [CF]).
 \bigskip
            
\noindent {\large\bf {Acknowledgements}}:
 I would like to thank D.~Sternheimer and Y.~Soibelman for their 
remarks and questions.

\bigskip

\noindent{\bf Mathematics Subject Classifications (1991)}:
{16S80, 55P35, 14F40, 14F25, 58A50}\\
{\bf Keywords:} {Deformation quantization, Hochschild cohomology,
little discs operads, 

\hskip16mm motives, quantum field theory.}

\begin{thebibliography}{BFFLS}
\bibitem[BFFLS]{BFFLS} Bayen, F, Flato, M., Fronsdal, C., Lichnerowicz, A., 
and Sternheimer, D.,
 {\sl Deformation theory and quantization. I. Deformations of symplectic
  structures}, Ann.\ Physics {\bf 111} (1978), no. 1, 61 - 110.

\bibitem[BV]{BV} Boardmann, J.~M., Vogt, R.~M., 
{\sl Homotopy invariant algebraic structures
   on topological spaces}, Springer-Verlag, Berlin, 1973,
    Lect. Notes in Math., Vol. 347.

\bibitem[BK]{BK} Broadhurst,D.~J. and  Kreimer, D.,  {\sl Association 
of multiple zeta values with positive knots via Feynman diagrams up 
to 9 loops},  {\tt hep-th/9609128}.

\bibitem[CF]{CF} Cattaneo, A. and  Felder,G.,  
{\sl A path integral approach to the Kontsevich 
quantization formula}, {\tt math/9902090}.

\bibitem[C]{C} Cohen, F.~R., {\sl The homology of $C_{n+1}$-spaces, $n\ge 0$},
 The homology of iterated loop spaces, Springer-Verlag, Berlin, 1976,
  Lect. Notes in Math., vol. 533, 207 - 351.

  
\bibitem[CK]{CK} Connes, A., Kreimer, D., 
{\sl Lessons from Quantum Field Theory}, Lett. Math. Phys. {\bf 48} (1999),
 this issue. 

\bibitem[DGMS]{DGMS} Deligne, P., Griffiths, Ph., Morgan, J., and  
Sullivan, D., {\sl Real homotopy theory of K\"ahler manifolds}, Invent. Math 
{\bf 29} (1975), 245 - 274.

\bibitem[D]{D} Drinfeld, V.~G., {\sl On quasi-triangular Quasi-Hopf algebras
 and a group closely related with $Gal({\overline {\Q}}/\Q)$},
   Leningrad Math. J., {\bf 2} (1991), 829 - 860.

\bibitem[EK]{EK} Etingof, P. and Kazhdan, D., {\sl Quantization of Lie 
Bialgebras, I},  Selecta Math., New Series, {\bf 2} (1996), no. 1, 1 - 41.
 
 \bibitem[FFS]{FFS} Flato, M., Fr\o nsdal, C., Sternheimer, D.,
  {\sl Singletons, physics in AdS universe and oscillations of 
composite neutrinos}, Lett. Math. Phys. {\bf 48} (1999),   this issue.

\bibitem[FM]{FM} Fulton, W. and MacPherson, R., {\sl Compactification of 
configuration spaces}, Ann.\ Math. {\bf 139} (1994), 183 - 225.

\bibitem[GV]{GV} Gerstenhaber, M. and Voronov, A.,
{\sl Homotopy $G$-algebras and moduli space
    operad}, Intern.\ Math.\ Res.\ Notices (1995), No. 3, 141 - 153.

\bibitem[GJ]{GJ} Getzler, E. and Jones, J.~D.~S., {\sl
  Operads, homotopy algebra and iterated
         integrals for double loop spaces},  {\tt hep-th/9403055}.

\bibitem[GK]{GK} Ginzburg, V. and Kapranov, M.,  {\sl Koszul duality for 
operads},  Duke Math. J. {\bf 76} (1994), no. 1, 203 - 272.

\bibitem[K1]{K1} Kontsevich, M., {\sl Deformation quantization of Poisson
 manifolds, I}, {\tt math/9709180}.

\bibitem[K2]{K2} Kontsevich, M., {\sl Feynman diagrams and low-dimensional 
topology}, First European Congress of Mathematics (Paris, 1992), Vol. II,
     Progress in Mathematics 120, Birkh\"auser (1994), 97 - 121.
     
 \bibitem[ML]{ML} Mac Lane, S.,
 {\sl Categories for working mathematician}, Springer-Verlag,
 Berlin, 1988.    

\bibitem[M]{M} May, J.~P., {\sl Infinite loop space theory}, Bull. Amer. 
Math. Soc. {\bf 83}  (1977), no. 4, 456 - 494.

\bibitem[MS]{MS} McClure, J. and  Smith, J.,  
{\sl Little 2-cubes and Hochschild cohomology}, preliminary announcement. 

\bibitem[Ne]{Ne} Nekov\'a\v r, J., {\sl Beilinson's Conjectures}, in
``Motives" (Editors U.~Jannsen, S.~Kleiman, J.~P.~Serre),
   Proceedings of Symposia in Pure Mathematics, AMS, vol. 55, part 1,
    537 - 570 (1994).


\bibitem[N]{N} Nori, M., private communication.

\bibitem[Q]{Q} Quillen, D., {\sl Homotopical algebra}, Lect. Notes Math., 
vol. 43, Berlin - Heidelberg - New York, Springer 1967.

\bibitem[S]{S} Serre, J.-P., 
{\sl Propri\'et\'es conjecturales des groups de Galois
 motiviques et des repr\'esentations $l$-adiques}, in
  ``Motives" (Editors U.~Jannsen, S.~Kleiman, J.-P.~Serre),
   Proceedings of Symposia in Pure Mathematics, AMS, vol. 55, part 1,
    377 - 400 (1994).

\bibitem[T1]{T1} Tamarkin, D.,  {\sl Another proof of M.~Kontsevich 
formality theorem}, {\tt math/9803025}.

\bibitem[T2]{T2} Tamarkin, D., {\sl Formality of Chain Operad of Small 
Squares},  {\tt math/9809164}.

\bibitem[V]{V} Voronov, A., {\sl The Swiss-Cheese Operad}, {\tt math/9807037}.

\bibitem[Z]{Z} Zagier, D., {\sl  Values of Zeta Functions and Their 
Applications}, 
First European Congress of Mathematics (Paris, 1992), Vol. II,
     Progress in Mathematics 120, Birkh\"auser (1994), 497 - 512.
\end{thebibliography}
\end{document}